\documentclass[12pt]{amsart}
\usepackage{amssymb, amsthm, amsfonts}
\usepackage{latexsym}
\usepackage{amsmath}
\usepackage{fullpage}
\usepackage{verbatim}
\usepackage{pslatex}

\input {cyracc.def}
\font\ququ=cmr10 scaled \magstep1
\font\tencyr=wncyr10 scaled \magstephalf
\font\tencyi=wncyi10 scaled\magstephalf
\font\tencysc=wncysc10 scaled\magstephalf
\def\rus{\tencyr\cyracc}
\def\rusi{\tencyi\cyracc}
\def\rusc{\tencysc\cyracc}

\newtheorem{thm}{Theorem}[section]
\newtheorem{lm}[thm]{Lemma}
\newtheorem{cl}[thm]{Corollary}
\newtheorem{prop}[thm]{Proposition}
\newtheorem*{probl}{Problem}

\theoremstyle{remark}
\newtheorem{rem}[thm]{Remark}
\newtheorem*{que}{Question}
\newtheorem*{prob}{Problem}

\theoremstyle{definition}
\newtheorem{opr}{Definition}
\newtheorem{ex}{Example}[section]

\newcommand{\re}[1]{\textrm {(\ref{#1})}}

\renewcommand{\theequation}{\thesection .\arabic{equation}}

\catcode`\@=\active
\catcode`\@=11
\def\@eqnnum{\hbox to
.01pt{}\rlap{\hskip-\displaywidth(\mathbf{\theequation})}}
\catcode`\@=12

\renewcommand{\iff}{if and only if }

\newcommand{\gt}{\mathfrak}

\newcommand {\ah}{\mathfrak{a}}

\newcommand {\g}{\mathfrak{g}}
\newcommand {\h}{\mathfrak{h}}

\newcommand {\el}{\mathfrak{l}}

\newcommand {\fN}{\mathfrak{N}}

\newcommand {\p}{\mathfrak{p}}
\newcommand {\q}{\mathfrak{q}}

\newcommand {\te}{\mathfrak{t}}

\newcommand {\z}{\mathfrak{z}}

\newcommand {\gln}{{\mathfrak{gl}_n}}
\newcommand {\sln}{{\mathfrak{sl}_n}}

\newcommand {\glv}{{\mathfrak{gl}(V)}}
\newcommand {\slv}{{\mathfrak{sl}(V)}}
\newcommand {\spv}{{\mathfrak{sp}(V)}}

\newcommand {\sov}{{\mathfrak{so}(V)}}

\newcommand {\esi}{\varepsilon}
\newcommand {\ap}{\alpha}

\newcommand {\vp}{\varphi}

\newcommand {\md}{/\!\!/}

\newcommand {\ad}{{\mathrm{ad\,}}}

\newcommand {\codim}{{\mathrm{codim\,}}}

\newcommand {\ind}{{\mathrm{ind\,}}}

\newcommand {\Ker}{{\mathrm{Ker\,}}}
\newcommand {\Ima}{{\mathrm{Im\,}}}

\newcommand {\rk}{{\mathrm{rk\,}}}
\newcommand {\spe}{{\mathrm{Spec\,}}}

\newcommand {\trdeg}{{\mathrm{trdeg\,}}}

\newcommand {\tri}{{\frak sl}_2}
\newcommand {\GR}[2]{{\textrm{{\bf #1}}}_{#2}}

\newcommand {\ov}{\overline}
\newcommand {\un}{\underline}

\newcommand {\rar}{\rightarrow}

\newcommand {\beq}{\begin{equation}}
\newcommand {\eeq}{\end{equation}}

\renewcommand{\le}{\leqslant}
\renewcommand{\ge}{\geqslant}

\newcommand{\zge}{\g_e}

\newfam\Bbbfam\newfam\eufam\newfam\eusfam%
\font\Bbbfont=msbm10 scaled 1200%
\font\frak=eufm10 scaled 1400%
\font\Bbbsmallfont=msbm8%
\textfont\Bbbfam=\Bbbfont\scriptfont\Bbbfam=\Bbbsmallfont%%
\font\euzw=eufm10 scaled 1200%
\font\euac=eufm7 scaled 1200%
\font\euacc=eufm7 scaled 1000%
\textfont\eufam=\euzw\scriptfont\eufam=\euac%
\scriptscriptfont\eufam=\euacc%
\font\euszw=eusm10 scaled 1200%
\font\eusac=eusm7 scaled 1200%
\font\eusacc=eusm7 scaled 1000%
\textfont\eusfam=\euszw\scriptfont\eusfam=\eusac%
\scriptscriptfont\eusfam=\eusacc%
\def\frak{\fam\eufam}%

\def\bbk{\hbox {\Bbbfont\char'174}}

\begin{document}
\setlength{\parskip}{2pt plus 4pt minus 0pt}
\hfill {\scriptsize March 1, 2005}
\vskip1ex
\vskip1ex

\title[The index of representations]{The index of representations
associated with stabilisers}
\author[D.\,Panyushev]{\sc Dmitri I. Panyushev}
\address[D.I.\,Panyushev]{ Independent University of Moscow,
Bol'shoi Vlasevskii per. 11, \
119002 Moscow, Russia}
\email{panyush@mccme.ru}
\author[ \ O.\,Yakimova]{Oksana S. Yakimova}
\address[O.S.\,Yakimova]{Max-Planck-Institut f\"ur Mathematik,
Vivatsgasse 7, \ 53111 Bonn, Germany}
\email{yakimova@mpim-bonn.mpg.de}
\thanks{This research was supported in part by
CRDF Grant no. RM1-2543-MO-03}
\maketitle
%\begin{center}
%{\footnotesize
%{\it Independent University of Moscow,
%Bol'shoi Vlasevskii per. 11 \\
%119002 Moscow, \quad Russia \\ e-mail}: {\tt panyush@mccme.ru }\\
%}
%\end{center}

\noindent
\section*{Introduction}

\noindent
The ground field $\bbk$ is algebraically closed and of characteristic zero.
For any finite-dimensional representation $\rho: \q\to \glv$ of a  Lie algebra
$\q$, one can define a non-negative integer which is called the {\it index\/}
of (the $\q$-module) $V$. Namely, if $V^*$ is the dual $\q$-module, then
\[
\ind (\q, V) = \dim V - \max_{\xi\in V^*}(\dim \q{\cdot}\xi)
\]
Here $\q{\cdot}\xi=\{ s{\cdot}\xi\mid s\in\q\}$ and
$s{\cdot}\xi$ is a shorthand  for $\rho^*(s)\xi$.
This definition goes back to {\sc Ra\"\i s}~\cite{rais}. 
Let $\q_v$ denote the stationary
subalgebra of $v\in V$.
For any $v\in V$, we can form the $\q_v$-module $V/\q{\cdot}v$. 
It was noticed by {\sc Vinberg} that one always have the inequality
\begin{equation}   \label{nervo}
   \ind(\q,V^*) \le \ind(\q_v, (V/\q{\cdot}v)^*) \ .
\end{equation}
The goal of this paper is to study conditions that guarantee us
the equality. If $V$ is the coadjoint representations of $\q$, then
the above index is equal to the index of $\q$ in the sense of Dixmier.
Here Vinberg's inequality reads
\[
   \ind\q \le \ind \q_\xi \quad \text{ for any }\ \xi\in \q^* \ .
\]
It is not always true that $\ind\q = \ind \q_\xi$, see Example~\ref{sl4} below.
However, it was conjectured by {\sc Elashvili} that if $\q=\g$ is semisimple, 
then this equality always holds. It easily seen that 
it suffices to prove the
equality $\ind\g = \ind\g_\xi$ only for the nilpotent elements 
$\xi\in \g\simeq \g^*$. The conjecture was recently proved by
{\sc Charbonnel}~\cite{char}. A proof for the classical Lie algebras,
with weaker assumptions on the ground field, 
was found independently by the second author~\cite{ksana1}.

One can consider two types of problems connected with Eq.~\re{nervo}.
{\it First\/}, to find properties of $v$ that guarantee the equality of the indices.
{\it Second\/}, to describe representations such that \re{nervo} turns into
equality for each $v\in V$.

We begin with pointing out two simple sufficient conditions.
If either $\q_v$ is reductive or $\dim\q_v{\cdot}v$ is maximal, then
Eq.~\re{nervo} turns into equality. 
Let $Q$ be a connected algebraic group with Lie algebra $\q$.
Given a representation $\rho: Q\to GL(V)$ (or $(Q:V)$ for short), we say that 
$(Q:V)$ has {\it good index behaviour\/} (GIB), if 
$ \ind(\q,V^*) =\ind(\q_v, (V/\q{\cdot}v)^*)$
 for each $v\in V$. We prove 
%%(surprisingly for us!) 
that most of sufficiently large reducible
representations have GIB. Namely, if $V$ is any (finite-dimensional
rational) $Q$-module, then $mV$ has GIB for any $m\ge \dim V$.
Another result of this sort asserts that if $V$ is a $Q$-module having GIB
and there is $\xi\in V^*$ such that $\q_\xi=0$, then $V\oplus W$ has GIB for
any $Q$-module $W$. It is also easily seen that any representation of an
algebraic torus has GIB.

Then we restrict ourselves to the case of reductive Lie algebras.
Here one can use the rich machinery and various tools of Invariant Theory.  
Let $G$ be a connected reductive group with Lie algebra $\g$.
Given a representation $\rho: G\to GL(V)$ (or $(G:V)$ for short), we say that 
$(G:V)$ has {\it good nilpotent index behaviour\/} (GNIB), if the equality
$ \ind(\g,V^*) =\ind(\g_v, (V/\g{\cdot}v)^*)$ holds
for any \un{nil}p\un{otent} element $v\in V$.
Using Luna's slice theorem, we prove that GIB is equivalent to that GNIB holds 
for any slice representation
of $(G:V)$. Furthermore, we prove that if $(G:V)$ is {\it observable\/}
(i.e., the number of nilpotent orbits is finite), then GNIB implies GIB.
As is well-known, the adjoint representation of $G$ is observable.

A related class of representations, with nice invariant-theoretic properties,
consists of the isotropy representations of symmetric pairs.
Since these representations are observable, it suffices to consider
the property of having GNIB for them. Let $(G,G_0)$ be a symmetric pair
with the associated $\mathbb Z_2$-grading $\g=\g_0\oplus\g_1$
and the isotropy representation $(G_0 : \g_1)$.
Abusing notation, we will say that $(G,G_0)$ has GNIB whenever the isotropy
representation has. A down-to-earth description of GNIB in the context of isotropy
representations is as follows.
Let $e\in\g_1$ be a nilpotent element, and $\g_e=\g_{e,0}\oplus\g_{e,1}$. 
Then GNIB property for $e$ means that the codimension of generic
$G_{e,0}$-orbits in $(\g_{e,1})^*$ equals the codimension of generic $G_0$-orbits
in $\g_1$, that is, the rank of the symmetric variety $G/G_0$. 
(By Vinberg's inequality, the first codimension cannot be less than the second one.)
It turns out that the analogue of the Elashvili conjecture (=\,Charbonnel's theorem)
does not always holds here, so that it is of interest to explicitly describe the
isotropy representations having GNIB.

In Sections~\ref{spv}-\ref{one-more-2}, we prove, using explicit matrix models, 
that the symmetric pairs
$(SL_n, SO_n)$, $(SL_{2n}, Sp_{2n})$, $(Sp_{2n}, GL_n)$, and $(SO_{2n}, GL_n)$
have GNIB. It is also shown that each symmetric pair of rank 1 has
GNIB, see Section~\ref{optimistic_end}.
On the other hand, we present a method of constructing isotropy 
representations without  GNIB, which makes use of even nilpotent orbits
of height 4. Combining this method with
the slice method, we are able to prove that most of the 
remaining isotropy representations do not have GNIB, see Section~\ref{ohne}.
As a result of our analysis and explicit 
calculations for small rank cases, we get a complete answer for the isotropy 
representations related to the classical simple Lie algebras.
The answer for $\sln$ is given below. 
%%$\mathbb Z_2$-gradings of $\sln$ (or $\gln$). 
\begin{thm}
Let $(SL_n,G_0)$ be a symmetric pair. Then it has GNIB if and only if
$\g_0$ belong to the following list:
(i) $\gt{so}_n$, (ii) $\gt{sp}_{2m}$ for $n=2m$, 
(iii) $\gt{sl}_m\times \gt{sl}_{n-m}\times\te_1$ with $m=1,2$, 
(iv) $\gt{sl}_3\times\gt{sl}_3\times\te_1$ for $n=6$.
[Here $\te_1$ stands for the Lie algebra of a one-dimensional torus.]
\end{thm}
%
%\\[.6ex]
\noindent
{\bf Acknowledgements.} {\ququ This paper was written during our stay at the
Max-Planck-Institut f\"ur Mathematik (Bonn). 
We are grateful to this institution for the warm hospitality and support.}

                %%%%
%%%%%%%%%%%%%%%%
%%%%%%%%%%%%%%%%   Section 1
%%%%%%%%%%%%%%%%
                %%%%

\section{The index of a representation}
\label{index}
\setcounter{equation}{0}

\noindent
Let $\q$ be a Lie algebra and
$\rho:\q \rar \glv$ a
finite-dimensional representation of $\q$, i.e., $V$ is a $\q$-module.
Abusing notation, we write $s{\cdot}v$ in place of $\rho(s)v$, if $s\in\q$
and $v\in V$. An element $v\in V$
is called {\it regular\/} or $\q$-{\it regular\/} whenever its
stationary subalgebra $\q_v=\{s\in\q\mid s{\cdot}v=0\}$ has minimal dimension. Because
the function  $v \mapsto \dim{\q}_v$ ($v\in V$) is upper semicontinuous,
the set of all $\q$-regular elements is open and dense in $V$.
%%Consider also the dual $\q$-module $V^\ast$.
%%\\[1ex]
\begin{opr}
The nonnegative integer
\[
\dim V - \max_{\xi\in V^*}(\dim \q{\cdot}\xi)=\dim V-\dim\q+\min_{\xi\in V^*}(\dim \q_\xi)
\]
is called the {\it index\/} of (the $\q$-module) $V$.
It will be denoted by $\ind (\q, V)$.
\end{opr}
Notice that
in order to define the index of $V$ we used elements of the dual
$\q$-module $V^*$.
This really makes a difference, since $\ind(\q,V)$ is not necessarily equal to
$\ind(\q,V^*)$ unless $\q$ is reductive.
\\[.7ex]
In case $\q$ is an algebraic Lie algebra, a more geometric description
is available. Let $Q$ be an algebraic group with Lie algebra $\q$.
Then $\ind\q=\dim\q-\max_{\xi\in\q^*}\dim Q{\cdot}\xi$.
By the Rosenlicht theorem \cite{R}, this number is also equal to
$\trdeg \bbk(V^*)^Q$.
Below, we always assume that $\q$ is algebraic, and consider $Q$ whenever it is
convenient.

If $v\in V$, then $\q{\cdot}v$ is a $\q_v$-submodule of $V$.
Geometrically, it is the tangent space of the orbit $Q{\cdot}v$ at $v$.
Then $V_v:=V/\q{\cdot}v$ is a $\q_v$-module as well.
By Vinberg's Lemma (see \cite[1.6]{cambr}),
we have
\begin{equation}   \label{vinb1}
\displaystyle \max_{x\in V}\dim(Q{\cdot}x)\ge
         \max_{\eta\in V_v}\dim(Q_v{\cdot}\eta) + \dim(Q{\cdot}v)
\end{equation}
for any $v\in V$.
It can be rewritten in equivalent forms:
\begin{gather}  
\label{vinb2}
     \trdeg \bbk(V)^Q \le \trdeg \bbk(V/\q{\cdot}v)^{Q_v} \quad \text{ or} \\
\label{vinb3}    
     \min_{x\in V}\dim(Q_x) \le \min_{\eta\in V_v}\dim((Q_v)_\eta)\quad \text{ or}   \\
\label{vinb4}
       \ind(\q,V^*) \le \ind(\q_v, (V/\q{\cdot}v)^*) \ .
\end{gather}
It is then natural to look for conditions that guarantee us the equality
%%in Eq.~\re{vinb4}.
This article is devoted to several aspects of the following problem
\begin{probl}
When does the equality hold in Eq.~\re{vinb1}--\re{vinb4} ?
\end{probl}
Every Lie algebra has a distinguished representation, namely,
the adjoint one. The index of the adjoint representation of
$\q$ is called simply {\it the index of\/} $\q$, denoted $\ind\q$.
That is, $\ind(\q,\q)=\ind\q$.
Let us take $V=\q^*$. Then $\q^*/\q{\cdot}\xi\simeq (\q_\xi)^*$ for any
$\xi\in \q^*$. Therefore
inequality \re{vinb4} in this situation reads
\begin{equation}   \label{vinb5}
     \ind\q \le \ind\q_\xi \quad \text{ for any }\ \xi\in \q^* .
\end{equation}
The coadjoint representation has some interesting
features. For instance, the $Q$-orbits in $\q^*$ are symplectic manifolds.
Hence $\ind\q_\xi-\ind\q$ is even for any  $\xi\in\q^*$.
However, even in this situation the inequality \re{vinb5} and
hence \re{vinb4} can be strict.
\begin{ex}    \label{sl4}
Let $\q$ be a  Borel subalgebra of $\gt{gl}_4$. It is well known
that $\ind \q=2$, see e.g. \cite[4.9]{mmj}. But there is a point
$\xi\in\q^*$ such that $\q_\xi$ is a 4-dimensional commutative
subalgebra, i.e., $\ind\q_\xi=4$. If $\q$ is represented as the space
of all upper-triangular matrices, then $\q^*\simeq \gt{gl}_4/[\q,\q]$
can be identified with
the space of all lower-triangular matrices. Then we take $\xi$ to be the
following matrix
$\xi=\begin{pmatrix} 0 &  \\
                     1 & 0 & \\
                     0 & 1 & 0 \\
                     0 & 0 & 1 & 0   \end{pmatrix}$.
\end{ex}
\noindent
Since the equality in Eq.~\re{vinb1}--\re{vinb4} does not always holds,
one has to impose some constraints  on $V$ and $Q$.
We begin with the following simple assertion, which is well known
to the experts.
\begin{prop}       \label{folk}
Suppose $Q_v$ is reductive. Then $\ind(\q,V^*) = \ind(\q_v, (V/\q{\cdot}v)^*)$.
\end{prop}\begin{proof}
In this case the $Q_v$-module $V$ is completely reducible, so that there
is a $Q_v$-stable complement of $\q{\cdot}v$, say $N_v$.
Let us form the associated fibre bundle
$Z_v:=Q\ast_{Q_v}N_v$. Recall that it is the (geometric) quotient of
$Q\times N_v$ by the $Q_v$-action defined by
$Q_v\times Q\times N_v\to Q\times N_v$, $(s,q,n)\mapsto (qs^{-1},s{\cdot}n)$.
The image of $(q,n)\in Q\times N_v$ in $Z_v$ is denoted by $q\ast n$.
Consider the natural $Q$-equivariant
morphism $\psi: Z_v\to V$, $\psi(q\ast n)=q{\cdot}(v+n)$. 
By construction, $\psi$ is \'etale in
$e\ast v\in Z_v$. It follows that the maximal dimensions of
$Q$-orbits in $V$ and $Z_v$ are the same, i.e.,
$\trdeg\bbk(Z_v)^Q=\trdeg\bbk(V)^Q$.
It remains to observe that
\[
    \max_{z\in Z_v}\dim(Q{\cdot}z) =
         \max_{\eta\in N_v}\dim(Q_v{\cdot}\eta) + \dim(Q{\cdot}v) \ ,
\]
which is a standard property of associated fibre bundles.
%%$\trdeg\bbk(Z_v)^Q=\trdeg\bbk(N_v)^{Q_v}$.
\end{proof}

\noindent
For the sake of completeness, we mention
the following obvious consequence of {\re{vinb1}}.

\begin{prop}  \label{max_orb}
If the dimension of $Q{\cdot}v$ is maximal, then
the action $(Q_v: V/\q{\cdot}v)$ is trivial and
the equality holds in \re{vinb1}.
\end{prop}

\begin{opr}   \label{GIB:def}
We say that the representation $(Q:V)$ has {\it good
index behaviour} ({\it GIB}, for short), if the equality 
\begin{equation}  \label{vinb=}
       \ind(\q,V^*) = \ind(\q_v, (V/\q{\cdot}v)^*) \
\end{equation}
holds for every $v\in V$. That is, inequality \re{vinb1} or \re{vinb4}
always turns into equality.
Another way is to say that $(Q:V)$ has GIB if and 
only if the function $v\mapsto \trdeg \bbk(V/\q{\cdot}v)^{Q_v}=
\ind(\q_v, (V/\q{\cdot}v)^*)$ is constant
on $V$.
\end{opr}%
\noindent
As an immediate consequence of Proposition~\ref{folk}, we obtain
\begin{prop}       \label{torus}
Let $Q$ be an algebraic torus. Then any $Q$-module has GIB.
\end{prop}
\noindent
For an arbitrary $Q$, it is not easy to prove that $V$ has
(or has not) GIB.
However, sufficiently ``large" reducible representations
always have GIB.
%%
%%Suppose that the representation $G:V$ is large, i.e., $\dim V\gg\dim\gt g$
%%and $\ind(\gt g, V^*)=\dim V-\dim\gt g$.
%%Suppose that for each $v\in V$ there is a generic $v_0$ such that 
%%$\gt gv\cap \gt v_0=\{0\}$. Then $G:V$ has GIB. We conjecture that this is 
%%always the case for ``very large'' representations.
%%Below we give two encouraging examples. 
%%
\begin{thm}        \label{big}
Let $\rho: Q\to GL(V)$ be an arbitrary linear  representation and $\dim V=n$.
Then the representation $(Q: mV^*)$ has GIB for any $m\ge n$. In this case,
$\ind(\q, mV^*)=nm-\dim\q$.
%%such that $\ind(\gt g, V^*)=\dim V-\dim\gt g$. 
%%Then $G:V^{(\dim V+1)}$ has GIB.
\end{thm}
\begin{proof} Our plan is to prove first the assertions
for $\q=\glv$, and then deduce from this the general case.

1) Assume that $\q=\glv$. It is clear that the generic stabiliser
for $(\glv: mV)$ is trivial for $m\ge n$, whence the equality for the index.

Let $\tilde v=(v_1,\ldots,v_m)$ be an arbitrary element of $mV$.
The rank of $\tilde v$, denoted $\rk\tilde v$, is the dimension of the 
linear span of the components $v_i$. If $\rk\tilde v=r\le n$, then without 
loss 
of generality one may assume that $\tilde v=(v_1,\ldots, v_r,0,\ldots,0)$,
where the vectors $v_1,\ldots, v_r$  form the part of the standard
basis for $V$. (Use the action of $GL_m$ that
permutes the coordinates of $\tilde v$.)
Then $GL(V)_{\tilde v}=\left(\begin{array}{cc}  I_r & \ast \\
                                         0   & \ast\end{array}\right)$
and $mV/\glv{\cdot}\tilde v\simeq (m-r)V$. It is easily seen that
$GL(V)_{\tilde v}$ has an orbit with trivial stabiliser in 
$mV/\glv{\cdot}\tilde v$.
This means that
$n(m-n)=\trdeg \bbk(mV)^{GL(V)}=
\trdeg\bbk(mV/\glv{\cdot}\tilde v)^{GL(V)_{\tilde v}}$ 
for any $\tilde v$, as required.

2) If $Q\subset GL(V)$ is arbitrary and  $\tilde v$ is as above, then 
$mV/\q{\cdot}\tilde v \supset (m-r)V$ and $\q_{\tilde v}\subset
\glv_{\tilde v}$. Hence $Q_{\tilde v}$ also has an orbit in
$mV/\q{\cdot}\tilde v$ with trivial stabiliser.
\end{proof}%
\begin{thm}   \label{GIB:sum}
Let $V$ be a $Q$-module having GIB
such that $\ind(\q, V^*)=\dim V-\dim\q$. 
Then for any $Q$-module $W$, $(Q:W\oplus V)$ has GIB and 
$\ind(\q, V^*\oplus W^*)=\dim V+\dim W-\dim\q$. 
\end{thm}
\begin{proof}
%%For each $x\in V$, its image in $V/\q{\cdot} v$ 
%%is denoted by $\overline x$. 
The assumption of having GIB and the equality for
$\ind(\q, V^*)$ mean that for any $v\in V$ there is $v_0$
such that $\dim (\q_v)_{\ov{v_0}}=0$,
where $\ov{v_0}$ stands for the image  of $v_0$ in $V/\q{\cdot}v$.
Our aim is to establish the similar property for $W\oplus V$.
Let $w+v\in W\oplus V$ be an arbitrary vector. 
%%Denote by $\overline{(w,v)}$ the image of 
%%$w+v$ in $(W\oplus V)/\gt g(w+v)$. 
Then 
\[
(\q_{v+w})_{\ov{(0,v_0)}}\subset
(\q_{v})_{\ov{(0,v_0)}}\subset
(\q_{v})_{\ov{v_0}}=\{0\},
\]
where $\ov{(0,v_0)}$ the image of $v_0$ in $(W\oplus V)/\q{\cdot}(w+v)$.
Therefore 
\[
\ind(\q_{v+w},(W\oplus V)/\q{\cdot}(w+v)^*)=
\dim (W\oplus V)-\dim\q{\cdot}(w+v)-\dim\q_{w+v}=
\dim (W\oplus V)-\dim\q. 
\]
Thus the function $v+w\mapsto\ind(\q_{v+w},(W\oplus V)/\q{\cdot}(w+v)^*)$ is 
constant, and we are done.
%%To conclude with, note that
%%$\ind(\q, (W\oplus V)^*)=\dim W+\dim V-\dim\q$.
\end{proof}%
\noindent
Combining the above theorems, we obtain

\begin{cl}     \label{GIB:big}
If\/ $V_1,V_2$ are arbitrary $Q$-modules and $m\ge \dim V_1$, then
$mV_1\oplus V_2$ has GIB.
\end{cl}

                %%%%
%%%%%%%%%%%%%%%%
%%%%%%%%%%%%%%%%   Section 2
%%%%%%%%%%%%%%%%
                %%%%

\section{Representations of reductive groups having GIB and GNIB }
\label{GIB:sect}
\setcounter{equation}{0}

\noindent
Let $G$ be a reductive algebraic group, and let $\rho: G\to GL(V)$ be
a finite-dimensional rational representation of $G$.
Recall that $v\in V$ is called {\it nilpotent\/}, if the closure of
the orbit $G{\cdot}v$ contains the origin, i.e.,
$\ov{G{\cdot}v}\ni 0$.
The set of all nilpotent elements is called the
{\it nullcone\/} and is denoted by $\fN(V)$.
Whenever we wish to stress that the nullcone depends on the group, we write $\fN_G(V)$.
A vector $v$ is said to be {\it semisimple}, if $\ov{G{\cdot}v}=G{\cdot}v$.
If $v$ is semisimple, then $G_v$ is reductive, and therefore the tangent
space $\g{\cdot}v \subset V$ has a $G_v$-stable complement, say $N_v$.
The natural representation $(G_v:N_v)$ is called the {\it slice representation\/}
(associated with $v$). We also say that it is a slice representation of
$(G{:}V)$. Notice that the initial representation itself can be regarded as
the slice representation associated with $0\in V$.
In this general situation, there is an analogue of the {\it Jordan decomposition\/},
which is well known for the elements of $\g$.
That is, for any $v\in V$ there are
a semisimple element $v_s$ and a (nilpotent) element $v_n$ such that
\begin{itemize}
\item \ $v=v_s+ v_n$;
\item \ $G_v \subset G_{v_s}$;
\item \ $v_n$ is nilpotent with respect to $G_{v_s}$, i.e., 
       $\ov{G_{v_s}{\cdot}v_n} \ni 0$.
\end{itemize}
This readily follows from Luna's slice theorem \cite{luna}.
Below, we recall how such a decomposition
is being constructed.
But, unlike the case of the adjoint representation,
a decomposition with the above properties is not unique. \\
As usual, $V\md G:=\spe \bbk[V]^G$ is the {\it categorical quotient\/} and
$\pi: V \to V\md G$ is the  {\it quotient mapping\/}.
Recall that $\fN_G(V)=\pi^{-1}\pi(0)$.

\begin{opr}   \label{GNIB:def}
We say that the representation $(G:V)$ has {\it good nilpotent
index behaviour} ({\it GNIB}, for short), if
equality \re{vinb=} holds
%%\begin{equation}  \label{vinbn=}
%%       \ind(\g,V^*) = \ind(\g_v, (V/\g{\cdot}v)^*) \
%%\end{equation}
for each {\sl nilpotent\/} element $v\in V$.
\end{opr}%
\noindent
First, we demonstrate that there are irreducible representations of reductive groups 
having no GNIB and thereby no GIB.
\begin{ex}   \label{contr}
Let $G=SL_2 \times SL_2$ and $V=R_3\otimes R_1$. Here $R_d$ stands for the simple
$SL_2$-module of dimension $d+1$. Hence $V$ is a simple $G$-module of dimension 8.
Let us show that $V$ has no GNIB.
A  generic stabiliser for this representation is finite, hence
$\ind(\g,V^*)=\dim V-\dim\g=2$.
As usual, we regard $R_d$ as the space of binary forms
of degree $d$. Let $(x^3, x^2y, xy^2, y^3)$ be a basis for $R_3$ and
$(u,z)$ a basis for $R_1$. Take $v=(x^3+ y^3)\otimes u$.
It is easily seen that $v$ is nilpotent. A direct computation shows that
the identity component of $G_v$ is 1-dimensional and unipotent.
However, the $\g_v$-module $V/\g{\cdot}v$  is trivial (and 3-dimensional).
Hence equality \re{vinb=} does not hold for $v$.
\end{ex}
\noindent
Our next goal is to understand a relationship between GIB and GNIB.
Clearly, if a representation has GIB, then it has GNIB as well.
As for the converse, we have the following general criterion.
%We do not know whether the converse is true. But the following weaker result
%holds.
%
\begin{thm}     \label{gnib+gib}
The representation $(G{:}V)$ has GIB if and only if 
every slice representation of $(G{:}V)$ has GNIB.
\end{thm}\begin{proof}
Actually, we prove a more precise statement.
Namely, suppose $v\in V$ is semisimple. Then 
equality \re{vinb=} is satisfied for every 
$y\in \pi^{-1}(\pi(v))$ if and only if 
the slice representation $(G_v{:}N_v)$ has GNIB.
\\[.7ex]
{\sf 1. "If" part.} \  
By Luna's slice theorem, $\pi^{-1}(\pi(v))\simeq G\ast_{G_v}\fN(N_v)$.
Therefore we may assume that $y=v+x$, where $x\in \fN(N_v)$.
This expression is just a Jordan decomposition for $y$, in the sense described above.
By assumption, we know that for any $x\in\fN(N_v)$ the following holds:
\begin{equation}   \label{for_x}
   \dim G_v{\cdot}x +\max_{\xi\in N_v/{\g_v}{\cdot}x}\dim (G_v)_x{\cdot}\xi=
   \max_{z\in N_v}\dim G_v{\cdot}z \ .
\end{equation}
Notice that $(G_v)_x=G_{v+x}=G_y$, since $y=v+x$ is a Jordan decomposition.
We want to show that
\begin{equation}  \label{for_y}
   \dim G{\cdot}y +\max_{\eta\in V/{\g}{\cdot}y}\dim (G_y){\cdot}\eta=
   \max_{z\in V}\dim G{\cdot}z \ .
\end{equation}
Again, since $y=v+x$ is a Jordan decomposition, we have $\dim G{\cdot}y=
\dim G{\cdot}v + \dim G_v{\cdot}x$. The following assertion is one of the many
consequences of Luna's slice theorem.

\begin{lm}    \label{isom}
The $G_y$-modules $N_v/{\g_v}{\cdot}x$ and $V/{\g}{\cdot}y$ are isomorphic.
\end{lm}\begin{proof}
First, we notice that both $N_v$ and $\g_v{\cdot}x$ are $G_y$-modules, since
$G_y=G_v\cap G_x$. Hence the first quotient is also a $G_y$-module.
Consider the $G$-equivariant morphism
\[
     %%\begin{array}{ccc}
         \psi:  G\ast_{G_v}N_v {\longrightarrow}V \ .
   %%              \downarrow & & \downarrow  \\
   %%  N_v\md G_v &  \overset{\psi\md G}{\longrightarrow} & V\md G   \end{array}
\]
%%where $\psi$ is , and the vertical arrows are quotient
%%mappings.
Recall that if $g{\ast} n\in G\ast_{G_v}N_v$ is an arbitrary point, then
$\psi(g{\ast} n):=g{\cdot}(v+n)$. Hence
$\psi(1{\ast} x) = y$. Set $\tilde y=1{\ast} x$.
It follows from the slice theorem that $G_{\tilde y}=G_y$
and
\[
T_{\tilde y}(G\ast_{G_v}N_v)/T_{\tilde y}(G{\cdot}\tilde y)
\to T_y V/T_v(G{\cdot}y)=V/{\g}{\cdot}y
\]
 is a $G_y$-equivariant bijection.
It remains to observe that the left-hand side is isomorphic to
$N_v/{\g_v}{\cdot}x$.
\end{proof}
\noindent
Thus, it follows from Lemma~\ref{isom} and the previous argument that the
left-hand side of \re{for_y} can be transformed  as follows
\begin{multline*}
  \dim G{\cdot}y +\max_{\eta\in V/{\g}{\cdot}y}\dim G_y{\cdot}\eta= \\
  =\dim G{\cdot}v + \bigl(\dim G_v{\cdot}x+
   \max_{\xi\in N_v/{\g_v}{\cdot}x}\dim G_y{\cdot}\xi\bigr)
   \overset{\re{for_x}}{=} \\
  \dim G{\cdot}v+\max_{z\in N_v}\dim G_v{\cdot}z \overset{Prop.\,\ref{folk}}{=}
\max_{z\in V}\dim G{\cdot}z \ ,
\end{multline*}
which completes the proof of the {\sf ``if"} part.
\\[.7ex]
{\sf 2. ``Only if" part.} \ Notice that the previous argument can be reversed.
\end{proof}
%
%There is a general philosophy of "good" properties for representations
%of reductive groups. 
In the light of the previous theorem, it is natural to ask the following 
natural 
\begin{que}
Is it true that "GNIB" implies "GIB" for any
representation of a reductive group?
\end{que}
\noindent 
We can give a partial answer to this question.
Recall that a representation $(G:V)$ is said to be {\it observable\/} if
the number of nilpotent orbits is finite.
This implies that each fibre of $\pi$ consists of finitely many orbits,
see e.g. \cite{kra}.
\begin{thm}    \label{oboz}
Suppose $(G:V)$ is observable. Then GNIB implies GIB.
\end{thm}\begin{proof}
Assume that this is not the case, i.e., $(G:V)$ has GNIB but
there is $v\in V$ such that $\ov{G{\cdot}v}\not\ni 0$
and Eq.~\re{vinb=} is not satisfied for $v$. 
We use the method of associated cones developed in \cite[\S\,3]{bokr}.
The variety $\ov{\bbk^* (G{\cdot}v)}\cap \fN(V)$ is the {\it associated cone\/} of
$G{\cdot}v$, denoted $\mathcal C(G{\cdot}v)$. It can be reducible, but each
irreducible component is of dimension $\dim G{\cdot}v$.
Let $G{\cdot}u$ be the orbit that is dense in an irreducible component of
$\mathcal C(G{\cdot}v)$. Here we use the hypothesis that $(G:V)$ is observable.
There is a 
morphism $\tau: \bbk \to \ov{\bbk^* (G{\cdot}v)}$
such that  $\tau(\bbk\setminus\{0\})\subset \bbk^* (G{\cdot}v)$ and 
$\tau(0)=u$.
%%sequence $\{x_i\}_{i=1}^\infty$ in $\bbk^* (G{\cdot}v)$ such that 
Since $\dim G{\cdot}u=\dim G{\cdot}v$, this implies that 
\begin{equation}   \label{limits}
\text{$\displaystyle\lim_{t\to 0}\g_{\tau(t)}=\g_u$ \ and \ 
$\displaystyle\lim_{t\to 0}\g{\cdot}\tau(t)=\g{\cdot}u$.}
\end{equation}
These two limits are taken in the suitable Grassmannians. 
By the assumption,
we have $\ind(\g_{\tau(t)}, (V/\g{\cdot}\tau(t))^*)< \ind (\g,V^*)$ for any $t\ne 0$.
In other words, 
\[
\max_{\eta\in V/\g{\cdot}\tau(t)} \dim G_{\tau(t)}{\cdot}\eta < 
\max_{z\in V} \dim G{\cdot}z -\dim G{\cdot}v \ .
\]
We claim that 
$\displaystyle \max_{\zeta\in V/\g{\cdot}u} \dim G_{u}{\cdot}\zeta \le
\max_{\eta\in V/\g{\cdot}\tau(t)} \dim G_{\tau(t)}{\cdot}\eta$.
%%=\max_{\eta\in V/\g{\cdot}x} \dim G_{x}{\cdot}\eta
This follows from the upper semi-continuity of dimensions of orbits and 
Eq.~\re{limits}.
%%the presence of the above limits.
The inequality obtained means that Eq.~\re{vinb=} is not satisfied for the nilpotent
element $u$. Hence $(G:V)$ has no GNIB, which contradicts the initial
assumption. This completes the proof.
\end{proof}%
\noindent
We do not know of whether the statement of Theorem~\ref{oboz}
remains true for arbitrary representations of $G$.

\noindent
Now, we turn to considering the adjoint representation of a reductive group $G$.
Here the condition of having GIB means that inequality~\re{vinb5} is,
in fact, equality. Because now $\g\simeq \g^*$, one may deal with centralisers
of elements in $\g$. As above, we write $\zge$ for the centraliser of $e\in\g$.
The following fundamental result was conjectured by Elashvili at the
end of 1980's and is recently proved by Charbonnel \cite{char}.

\begin{thm} (Charbonnel)    \label{alela:conj}
The adjoint representation of a reductive group $G$ has GNIB. In other words, if
$e\in\g$ is a nilpotent element, then $\ind \zge=\rk \g$.
\end{thm}%
\noindent
In \cite{ksana1}, this theorem is independently proved for the classical Lie 
algebras. Some partial results for "small" orbits were obtained earlier in
\cite{cambr} and~\cite{amazing}.

A remarkable fact is that, for the adjoint representation, each slice
representation is again the adjoint representation (of a centraliser).
Hence Theorems~\ref{gnib+gib} and \ref{alela:conj} readily imply that the
adjoint representation has GIB. Another way to deduce GIB is to refer to 
Theorems~\ref{oboz} and \ref{alela:conj}, and the fact that the adjoint 
representation is observable.
%%Thus, we obtain the first and, at the moment, unique series of
%%representations with GIB.
\\[.7ex]
From the invariant-theoretic point of view, adjoint representations
have the best possible properties. Isotropy representations of symmetric spaces
form a class with close properties. So, it is natural to inquire whether
%%properties of having GIB and GNIB for 
these representations have GIB and GNIB.
Recall the necessary setup.

Let $\sigma$ be an involutory automorphism of $\gt g$.
Then $\g=\g_0\oplus \g_1$ is the direct sum of $\sigma$-eigenspaces.
Here $\g_0$ is a reductive subalgebra and
$\g_1$ is a $\g_0$-module. Write $G_0$ for the connected subgroup of
$G$ with Lie algebra $\g_0$. With this notation, our object of
study is $(G_0{:}\g_1)$,
the {\it isotropy representation\/} of the symmetric pair $(G, G_0)$.
By \cite{KR}, these representations are observable, so that Theorem~\ref{oboz}
applies. Therefore we will not distinguish the properties GIB and GNIB in the context
of isotropy representations of symmetric pairs.
In the rest of the paper, we deal with the following

\begin{prob}
For which involutions $\sigma$ the representation $(G_0{:}\g_1)$ has GNIB ?
\end{prob}
\noindent
For future use, we record the following result.

\begin{lm} \label{z2}
Let $\q=\q_0\oplus\q_1$ be an arbitrary $\mathbb Z_2$-graded
Lie algebra and $\q^*=\q_0^*\oplus\q_1^*$ the corresponding decomposition
of the dual space. For any $\xi\in\q^*_1$ the stationary subalgebra
$\q_\xi$ posesses the induced $\mathbb Z_2$-grading and \ 
%%we have
$\dim\q_0-\dim\q_1=\dim\q_{\xi,0}-\dim\q_{\xi,1}$.
\end{lm}%
\begin{proof}
This claim is well known if $\q$ is reductive and one identifies
$\q$ and $\q^*$, see \cite[Prop.\,5]{KR}. 
The general proof is essentially the same.
\end{proof}%
\noindent
Let us give an interpretation of GNIB for the
isotropy representations, which is helpful in practical
applications. 
It is known that $x\in\g_1$ is nilpotent 
in the sense of the above definition (i.e., as an element of the $G_0$-module $\g_1$)
if and only if it is nilpotent as
an element of  $\g$.
Formally, $\fN_{G_0}(\g_1)=\fN_{G}(\g)\cap \g_1$. If $e\in \fN(\g_1)$,
and $\zge=\g_{e,0}\oplus\g_{e,1}$ is the induced ${\mathbb Z}_2$-grading,
%%%then $\zge$ inherits the ${\mathbb Z}_2$-grading:
%%%$\zge=\g_{e,0}\oplus\g_{e,1}$, and 
then $\g_{e,0}$ is precisely the stationary
subalgebra of $e$ in $\g_0$.
Now, inequality~\re{vinb4} reads
\[
   \ind(\g_0, (\g_1)^*)\le \ind (\g_{e,0}, (\g_1/[\g_0,e])^*) \ .
\]
Using a $G$-invariant inner product on $\g$,
one easily shows that $\g_1/[\g_0,e]\simeq (\g_{e,1})^*$.
Recall also that $\g_1$ is an orthogonal $G_0$-module, i.e.,
$G_0\to SO(\g_1)$. The number
$\ind(\g_0, (\g_1)^*)=\ind(\g_0, \g_1)$ equals the Krull dimension
of the invariant ring $\bbk[\g_1]^{G_0}$, which in turn
is equal to the rank of $G/G_0$ (in the sense of the theory of
symmetric varieties).
Thus, we obtain

\begin{prop}   \label{z2-gnib} $\phantom{svolochi!}$

1. For any $e\in \fN(\g_1)$, we have
$\rk(G/G_0) =\ind(\g_0,\g_1) \le \ind (\g_{e,0}, \g_{e,1})$.

2. The following conditions are equivalent:
\begin{itemize}
\item[\sf(i)] \ 
The isotropy representation $(G_0{:}\g_1)$ has GNIB;
\item[\sf (ii)] \ 
%%if and only if  
for any  $e\in \fN(\g_1)$ we have
$\rk(G/G_0)=\ind (\g_{e,0}, \g_{e,1})$;
\item[\sf (iii)] \ for any $e\in \fN(\g_1)$ there is an 
$\alpha\in\gt g_e^*$ such that $\alpha(\g_{e,0})=0$ and 
$\dim(\gt g_{e,1})_\alpha=\rk(G/G_0)$.
\end{itemize}
\end{prop}
\begin{proof} Part 1 and the equivalence of (i) and (ii)
follow from the previous discussion.
To prove the equivalence of (ii) and (iii), we note that
if $\alpha(\g_{e,0})=0$, then $\ap$ can be regarded as an element of
$(\g_{e,1})^*$. Then 
\[
\codim \g_{e,0}{\cdot}\ap=
\dim\gt g_{e,1}-\dim\gt g_{e,0}+\dim(\gt g_{e,0})_\alpha \overset{Lemma~\ref{z2}}{=} 
%%\dim(\gt g_{e,1})_\alpha-\dim(\gt g_{e,0})_\alpha+\dim(\gt g_{e,0})_\alpha=
\dim(\gt g_{e,1})_\alpha.
\]
Hence, $\ind(\gt g_{e,0},\gt g_{e,1})=\min\dim(\gt g_{e,1})_\alpha$, 
where minimum is taken over all $\alpha\in\gt g_e^*$ 
such that $\alpha(\g_{e,0})=0$.
\end{proof}

\noindent
Below, we show that there are isotropy representations with and without GNIB. 
%%On the other hand,
%%there are series of involutions which do provide representations with GIB.

                %%%%
%%%%%%%%%%%%%%%%
%%%%%%%%%%%%%%%%   Section 3
%%%%%%%%%%%%%%%%
                %%%%

\section{Isotropy representations for the outer involutions 
of $\glv$}
\label{spv}
\setcounter{equation}{0}

\noindent
Let $V$ be a finite-dimensional vector space over $\bbk$.
If $\sigma$ is an outer involution of $\slv$, then
$\g_0$ is isomorphic to either $\spv$ or $\sov$. Of course, the first case
is only possible if $\dim V$ is even. It will technically be easier to
deal with $\g=\glv$ and assume that the centre of $\glv$ lies in $\g_1$.
Then the $Sp(V)$-module
$\g_1$ is isomorphic to $\wedge^2 V$ and the $SO(V)$-module
$\gt g_1$ is isomorphic to ${\mathcal S}^2V$. The goal of
this section is to prove that the isotropy representations
$(Sp(V):\wedge^2 V)$ and $(SO(V):{\mathcal S}^2 V)$ have GIB.

Recall the necessary set-up.
Let $(\phantom{,},\phantom{,})$ be a non-degenerate
symmetric or skew-symmetric form on  $V$; that is,
$(v,w)=\esi(w,v)$, where $v,w\in V$ and $\esi=+1$ or $-1$.
Let $J$ denote the matrix of $(\,,\,)$ with respect to some basis of $V$.
Then $(v,w)=v^tJw$,
where $v,w$ are regarded as column vectors and
the symbol $(\phantom{,})^t$ stands for the transpose. 
Since $J^t=\pm J$, the mapping $A \to \sigma(A):=-J^{-1} A^t J$ is an involution of 
$\glv$.
Let $\glv=\g_0\oplus\gt g_1$ be the corresponding $\mathbb Z_2$-grading.
Here $\g_0$ consists of the linear transformations preserving the form $(\,,\,)$,
i.e., satisfying the property $(vA,w)=-(v,Aw)$ for all $v,w\in V$.
The elements of $\gt g_1$ multiply  the form 
$(\phantom{,},\phantom{,})$  by $-1$, i.e.,
\begin{equation}   \label{3.1}
(A v,w)=(v,A w)  \quad \text{ for all $A\in\gt g_1$ and $v,w\in V$.}
\end{equation}
Recall standard facts concerning nilpotent elements
in $\g=\glv$, mainly in order to fix the notation.
\\[.6ex]
Let $e\in\g$ be a nilpotent element and $m=\dim \Ker(e)$. 
By the theory of Jordan normal form,  there are vectors
$w_1,\ldots,w_m\in V$ and non-negative integers $d_1,\ldots,d_m$
such that $e^{d_i+1}{\cdot}w_i=0$
and $\{ e^s{\cdot}w_i \mid 1\le i\le m, \ 0\le s\le d_i\}$ is
a basis for $V$. 
Set $V_i=\langle w_i,e{\cdot}w_i,\ldots, e^{d_i}{\cdot}w_i\rangle$
and $W=\langle w_1,\ldots,w_m\rangle$. Then $V=\oplus_{i=1}^m V_i$
and $V=W\oplus \Ima(e)$.
The spaces $\{V_i\}$ are called the Jordan (or {\it cyclic}) spaces of the 
nilpotent element $e$.
%%Let $V_i\subset V$ be the space of the $i$-th Jordan block,
%%so that $V=\oplus_{i=1}^m V_i$. Set $d_i+1=\dim V_i$.
%%In each $V_i$, we choose a cyclic vector $w_i$; that is,
%%$\{w_i,e{\cdot}w_i,\ldots, e^{d_i}{\cdot}w_i\}$ is a basis for $V_i$.
%We shall say that $W$ is {\it a space of cyclic vectors}.

Suppose $\varphi\in\g_e$. Because
$\varphi(e^s{\cdot}w_i)=e^s{\cdot}\varphi(w_i)$, the linear map $\varphi$
is determined by its values on $W$.
In other words, if
\[ 
\varphi(w_i)=\sum_{j,s}
c_i^{j,s}(e^s{\cdot}w_j), \mbox{ where } c_i^{j,s}\in\bbk \ ,
\]
then $\varphi$ is determined by the coefficients
$c_i^{j,s}=c_i^{j,s}(\varphi)$. In what follows, we will only indicate the
values of $\vp$ on the cyclic vectors $\{w_i\}$.

A basis for $\g_e$ consists of the maps
$\{\xi_i^{j,s}\}$ given by
\[
\xi_i^{j,s}:\left\{
\begin{array}{llcl}
w_i & \mapsto & e^s{\cdot}w_j & \\
w_t & \mapsto & 0  & \mbox{ if } t\ne i \\
\end{array}\right. ,\quad \mbox{where} \quad  1\le i,j\le m \ \mbox{ and }\
\max\{d_j-d_i, 0\} \le s\le d_j \ .
\]

%We will be interested in the subalgebra $\gt h$ of all
%elements $\vp \in \gt z(e)$ preserving the Jordan blocks of
%$e$, i.e., having the property that $\vp(V_i)\subset V_i$
%for each $i$. Obviously, this happens \iff
%$c_i^{j,s}(\varphi)=0$ unless $i=j$. Another way to express
%this is to say that a basis for this space consists of all
%$\xi_i^{i,s}$, $0\le s\le d_i$.

\begin{lm}   \label{restr}
In the above setting, suppose that $e\in\fN(\g_1)$.
Then the cyclic vectors $\{w_i\}$ and thereby the spaces $\{V_i\}$ can be chosen
such that the following conditions are satisfied:
\begin{itemize}
\item[\sf (i)] \ If\/ $\esi=-1$, then 
the set $\{1,2,\ldots,m\}$ can be partitioned in pairs $(i,i^*)$ such that
$d_i=d_{i*}$ and
$w_i$ is orthogonal to all basis vectors $e^s{\cdot}w_j$ except
$e^{d_i}{\cdot}w_{i^*}$.
(\mbox{\rm Here $i\ne i^*$}.)
\item[\sf (ii)] \  If\/ $\esi=1$, then $(V_i, V_j)=0$ for $i\ne j$ and
the restriction of $(\,,\,)$ to each $V_i$ is non-degenerate.
\end{itemize}
\end{lm}%
\begin{proof} 
We argue by induction on $m=\dim\Ker(e)$.
\\[.6ex]
It follows from Eq.~\re{3.1} that
$\Ker(e^i)$ and $\Ima(e^i)$ are orthogonal with
respect to $(\,,\,)$. In particular, 
$\Ker(e)$ is orthogonal to $\Ima(e)$, and 
$(\,,\,)$ induces a
non-degenerate pairing between $W$ and $\Ker(e)$.
%%Since $e\in\g_1$, we have
%%$(w_i,e^s{\cdot}v)=(e^s{\cdot}w_i,v)$ and therefore  
Suppose $d_1=\min_i\{d_i\}$. There is a vector 
$e^{d_i}{\cdot} w_i\in\Ker(e)$ for some $i$
such that $(w_1,e^{d_i}w_i)\ne 0$.
Then $d_i\le d_1$, hence $d_i=d_1$ in view of the minimality of $d_1$.
\\
The rest of the argument splits.

(i) In the symplectic case ($\esi=-1$), we have
\[
  (w_1,e^{d_1}{\cdot}w_1)=(e^{d_1}{\cdot}w_1,w_1)
=-(w_1,e^{d_1}{\cdot}w_1)=0.
\]
Hence, $i\ne 1$. It is easily verified that the 
restriction of $(\,,\,)$ to either $V_1$ or $V_i$ is zero,
while the restriction to $V_1\oplus V_i$ is 
non-degenerate. Therefore, we may take $1^*=i$. Then 
all other $w_j$ can be chosen in $(V_1\oplus V_i)^{\perp}$,
the $e$-invariant orthogonal complement to $V_1\oplus V_i$. 
%%%this procedure yields a required set of cyclic vectors $\{w_i\}$.
%%such that $w_i$ is orthogonal to all vectors $e^sw_j$ except
%%$e^{d_i}w_{i^*}$, where $d_i=d_{i*}$, $i\ne i^*$ and $(i^*)^*=i$.

(ii) Consider the orthogonal case ($\esi=1$). If $i=1$, then the restriction of
$(\,,\,)$ to $V_1$ is non-degenerate and we may choose the
remaining cyclic vectors in $V_1^{\perp}$. %%\cap W$. 
If $i\ne 1$ and $(\,,\,)$
is degenerate on both $V_1$ and $V_i$, then we make the following modification
of $w_1$ and $w_i$. Our assumption implies that
$(w_1,e^{d_1}{\cdot}w_1)=0$ and $(w_i,e^{d_1}{\cdot}w_i)=0$.
Set $w_1':=w_1+ w_i$ and $w_i':=w_1-w_i$. Then
$(w_1',e^{d_1}{\cdot}w_1')=2(w_1,e^{d_1}{\cdot}w_i)\ne 0$ and
$(w_i',e^{d_1}{\cdot}w_i')=-2(w_1,e^{d_1}{\cdot}w_i)\ne 0$.
This means that the restriction of $(\,,\,)$ to the 
cyclic space generated by either $w'_1$ or $w'_i$ is non-degenerate.
%%\{w_1',\ldots, e^{d_1}{\cdot}w_1'\}$ is non-degenerate.
\end{proof}

\begin{thm}     \label{so-gib}
The representation $(SO(V) : {\mathcal S}^2 V)$ has GNIB.
\end{thm}
\begin{proof} Here $\rk(G/G_0)=\dim V$. Let $e\in \fN(\g_1)$.
Recall that  $\sigma$ induces the decomposition
$\g_e=\g_{e,0}\oplus\g_{e,1}$.
We choose the cyclic vectors for $e$ as described in Lemma~\ref{restr}(ii).
Define  $\alpha\in(\g_e)^*$ by
\[
\alpha(\varphi)=\sum_{i=1}^m a_i c_i^{i,d_i},
\]
where  $c_i^{j,s}$ are the coefficients of $\vp$ 
and $\{a_i\}$ are pairwise different non-zero numbers.
Then 
$(\g_e)_\ap$ consists of all maps  
in $\g_e$ preserving the Jordan spaces $V_i$ \ \cite[Sect.\,2]{ksana1}, i.e.,
\[
  (\g_e)_\ap=\langle \xi_i^{i,s}\mid 1\le i\le m, \ 0\le s\le d_i\rangle \ ,
\]
where $\langle\ldots\rangle$ denotes the $\bbk$-linear span. 
Hence $\dim(\g_e)_\ap =\sum_i (d_i+1)=
\dim V$. 
%%% and  $(\g_e)_\ap\subset \g_{e,1}$.
We claim that  $\ap(\g_{e,0})=0$. Indeed, assume the converse,
i.e., $\vp=\sum c_i^{j,s}\xi_i^{j,s}\in \g_{e,0}$ and $\ap(\vp)\ne 0$.
This means that $c_i^{i,d_i}\ne 0$ for some $i$.
Then 
\[
(\vp(w_i),w_i)=(w_i,\vp(w_i))=
c_i^{i,d_i}(w_i,e^{d_i}{\cdot}w_i)\ne 0,
\]
which in  view of Eq.\re{3.1} 
contradicts the fact that $\vp\in \g_{e,0}$.
Hence $\ap\in (\g_{e,1})^*\subset (\g_e)^*$.
By Lemma~\ref{z2}, we have
\[  
\dim(\g_{e,1})_\alpha-\dim(\g_{e,0})_\alpha=
\dim \g_{e,1}-\dim \g_{e,0}=\dim\gt g_1-\dim\gt g_0=\dim V \ .
\]
On the other hand, 
\[
\dim(\g_{e,1})_\alpha+\dim(\g_{e,0})_\alpha=\dim(\gt g_e)_\alpha=
\dim V.
\]
Hence $\dim(\g_{e,1})_\alpha=\dim V=\rk G/G_0$. Thus, $(SO(V) :{\mathcal S}^2 V)$ 
has GNIB in view of Proposition~\ref{z2-gnib}(2). 
\end{proof}

\begin{thm}    \label{sp-gib}
The representation $(Sp(V) : \wedge^2 V)$ has GNIB.
\end{thm}
\begin{proof} Put $\dim V=2n$. It is well known that 
$\rk(G/G_0)=\dim V/2=n$.
\\[.6ex]
Let $e\in \fN(\g_1)$.
We choose the cyclic vectors for $e$ as described in Lemma~\ref{restr}(i).
By Eq.~\re{3.1}, we have
$(w_i, e^{d_i}{\cdot}w_{i^*})=(e^{d_i}{\cdot}w_i, w_{i^*})=-
(w_{i^*}, e^{d_i}{\cdot}w_{i})$.
 This implies that 
$\xi_i^{i,d_i}+\xi_{i^*}^{i^*, d_i}
\in \g_{e,1}$ and $\xi_i^{i,d_i}-\xi_{i^*}^{i^*, d_i}
\in \g_{e,0}$ for each $i$.

Define  $\alpha\in(\g_e)^*$ by
\[
\alpha(\varphi)=\sum\limits_{i=1}^m a_i c_i^{i,d_i},
%%\enskip a_i\in\bbk,
\]
where  $c_i^{j,s}$ are the coefficients of $\vp$ and $\{a_i\}$ are non-zero
numbers such that $a_i=a_j$ \iff $i=j^*$. A direct computation shows that
$(\g_e)_\alpha$ consists of all elements of $\g_e$
preserving the subspaces $V_i\oplus V_{i^*}$ for each pair $(i,i^*)$, 
cf. \cite[Sect.\,2]{ksana1}. More concretely,
$$
 (\g_e)_\alpha=\langle \xi_i^{i,s}\mid 1\le i\le m, \ 0\le s\le d_i\rangle 
\oplus
  \langle \xi_i^{i^*,s}\mid 1\le i\le m, \ 0\le s\le d_i\rangle \ .
$$
Hence,
$\dim(\g_e)_\alpha=2\dim V=4n$. 
%%We denote the restriction of $\alpha$ to $\g_{e,1}$ by the same symbol.
As in the proof of Theorem~\ref{so-gib}, one can show that $\alpha\vert_{\g_{e,0}}=0$. 
Therefore $\ap$ can be regarded as an element
of  $(\g_{e,1})^*\subset (\g_e)^*$. 
Using Lemma~\ref{z2}, we obtain
$$ 
\dim(\g_{e,0})_\alpha-\dim(\g_{e,1})_\alpha=
\dim \g_{e,0}-\dim \g_{e,1}=\dim\gt g_0-\dim\gt g_1=2n \ .
$$
It follows that $\dim(\g_{e,0})_\ap=3n$.
and $\dim(\g_{e,1})_\ap=n$.
%%%Thus, %%%using the definition of the index and Proposition~\ref{z2-gnib},
By Proposition~\ref{z2-gnib}(2),
we conclude 
%%\[
%%\ind(\g_{e,0},\g_{e,1})\le \dim\g_{e,1}-\dim\g_{e,0}+\dim(\g_{e,0})_\ap=-2n+3n=n=
%%\rk G/G_0 \le \ind(\g_{e,0},\g_{e,1}) ,
%%%\]
%%which means 
that $(Sp(V) : \wedge^2 V)$ has GNIB. 
\end{proof}
\noindent
In Section~\ref{ohne}, we show that most of the isotropy representations associated 
with inner involutions of $\glv$ do not have GNIB.

                %%%%
%%%%%%%%%%%%%%%%
%%%%%%%%%%%%%%%%   Section 4
%%%%%%%%%%%%%%%%
                %%%%

\section{The isotropy representation of $(\gt{sp}_{2n},\gt{gl}_n)$}
%%%a symmetric pair associated with $\spv$}
\label{one-more}
\setcounter{equation}{0}

\noindent
In this section, $\dim V=2n$, $\gt g=\gt{sp}(V)$, and $(\,,\,)$ is a
$\g$-invariant skew-symmetric form on $V$.  
Let $\sigma$ be an involution of $\g$ such that $\g_0\simeq
\gln$. This can explicitly be described as follows. 
%%%and $\gt g_0=\gt g^\sigma=\gt{gl}_n$. 
Let $V=V_+\oplus V_-$ be a Lagrangian decomposition of $V$.
Then $G_0$ can be taken as the subgroup
of $G=Sp(V)$ preserving this decomposition. 
Here $G_0\simeq GL(V_+)$, $V_-\simeq (V_+)^*$ as $G_0$-module, and
the $G_0$-module $\g_1$ is isomorphic to 
${\mathcal S}^2 V_+\oplus ({\mathcal S}^2 V_+)^*$.

Keep the notation introduced in the previous section.
In particular, $V_i$, $i=1,\ldots,m$, are the Jordan spaces of $e\in\fN(\glv)$,
$\dim V_i=d_i+1$, and $w_i\in V_i$ is a cyclic vector.

\begin{lm}   \label{restr2}
Let $e\in\fN(\g_1)$.
Then the cyclic vectors $\{w_i\}_{i=1}^m$ and hence the $\{V_i\}$'s 
can be chosen such that the following properties are satisfied:
\begin{itemize}
\item[\sf (i)] \ 
there is an involution $i\mapsto i^*$ on the set
$\{1,\dots, m\}$ such that 
\begin{itemize}
\item[$\circ$] \ $d_i=d_{i^*}$, 
\item[$\circ$] \ $i=i^*$ if and only if $d_i$ is odd, 
\item[$\circ$] \ $(V_i, V_j)=0$ if $i\ne j,j^*$;
\end{itemize}
\item[\sf (ii)] \ $\sigma(w_i)=\pm w_i$.
\end{itemize}
\end{lm}
\noindent 
The proof 
%%, which is an exercise in Linear Algebra,
is left to the reader (cf. the proof of Lemma~\ref{restr}).
Actually, part (i) is a standard property of nilpotent orbits in $\spv$.
Then part (ii) says that in the presence of the involution $\sigma$
the cyclic vectors for $e\in\fN(\g_1)$ can be chosen to be $\sigma$-eigenvectors.

\begin{thm}   \label{more_gib}
The representation 
$(GL(V_+) : {\mathcal S}^2 V_+\oplus ({\mathcal S}^2 V_+)^*)$ 
has GNIB. 
\end{thm}
\begin{proof} 
Recall that $\gt{sp}(V)$ is a symmetric 
subalgebra of $\widetilde{\gt g}:=\glv$.  
Let $\widetilde{\gt g}=\gt{sp}(V)\oplus\widetilde{\gt g}_1$ be the corresponding
$\mathbb Z_2$-grading. Then we have a hierarchy of involutions:
\\[.6ex]
\centerline{
$\tilde\g=\tilde\g_0\oplus\tilde\g_1$ \ and \ $\tilde\g_0=\g=\g_0\oplus\g_1$.
}
\\[.6ex]
Let $e\in\fN(\gt g_1)$. In this case, we have $\rk (G/G_0)=n$.
Hence, by Proposition~\ref{z2-gnib} our goal is to find an element 
$\ap\in (\g_{e,1})^*$ such that
$\dim(\g_{e,1})_\ap=n$.
Let $\widetilde{\gt g}_e$ and $\widetilde{\gt g}_{e,1}$ 
denote the centraliser of $e$ in 
$\widetilde{\gt g}$ and $\widetilde{\gt g}_1$, respectively.
In view of the above hierarchy, we have 
\\[.6ex]
\centerline{
$\tilde\g_e=\g_{e}\oplus \tilde\g_{e,1}=
\g_{e,0} \oplus \g_{e,1}\oplus \tilde\g_{e,1}$.
}
\\[.6ex]
Choose the cyclic vectors for $e$ as prescribed by
Lemma~\ref{restr2}. 
We normalise these vectors such that 
$(w_i,e^{d_i}{\cdot}w_i)=1$ if $i=i^*$; and
$(w_i, e^{d_i}{\cdot}w_{i^*})=-(w_{i^*}, e^{d_i}{\cdot}w_i)=\pm 1$ if $i\ne i^*$. 
Then $\gt g_e$ has a basis 
$\xi_i^{j,d_j-s}+\esi(i,j,s)\xi_{j^*}^{i^*,d_i-s}$, 
where $\esi(i,j,s)=\pm 1$ depending on $i,j$ and $s$; and 
$\xi_i^{j,d_j-s}-\esi(i,j,s)\xi_{j^*}^{i^*,d_i-s}$ form a basis 
for $\widetilde{\gt g}_{e,1}$. 

Define  $\alpha\in (\tilde{\g}_e)^*$ by 
\[
\alpha(\varphi)=\left(\sum\limits_{i, \ i=i^*} a_ic_i^{i,d_i}\right)+
                \sum\limits_{(i,i^*),\ i\ne i^*}
 a_i(c_i^{i^*,d_i}+c_{i^*}^{i,d_i}), 
\]
where $c_i^{j,s}$ are coefficients of $\varphi\in\widetilde{\gt g}_e$,
$a_i=a_{i^*}$, and $a_i\ne\pm a_j$ if $i\ne j,j^*$.
The stationary subalgebra $(\widetilde{\gt g}_e)_\alpha$ consists of all
maps preserving cyclic spaces generated by $w_i$ for $i=i^*$
and $w_i+w_{i^*}$, $w_i-w_{i^*}$ for $i\ne i^*$.
More precisely,
\[
\begin{array}{l}
(\widetilde{\gt g}_e)_\alpha = 
 \langle \xi_i^{i,s}\mid 1\le i\le m, \ i=i^*, \ 0\le s\le d_i\rangle 
\oplus \\
\qquad \qquad \ \ \oplus
  \langle \xi_i^{i,s}+\xi_{i^*}^{i^*,s}, \ \xi_i^{i^*,s}+\xi_{i^*}^{i,s}
\mid 1\le i\le m, \ i\ne i^*, \ 0\le s\le d_i\rangle \ . \\
\end{array}
\]
First, we show that $\alpha(\widetilde{\gt g}_{e,1})=0$.
Assume that  
 $\alpha(\xi_i^{j,d_j-s}-\esi(i,j,s)\xi_{j^*}^{i^*,d_i-s})\ne 0$
for some $\xi_i^{j,d_j-s}-\esi(i,j,s)\xi_{j^*}^{i^*,d_i-s}\in\widetilde{\gt g}_{e,1}$.
Then $j=i^*$, $s=0$ and $\esi(i,i^*,0)=-1$. 
But $\xi_i^{i^*,d_i}\in\gt g$ for all $i$, 
%%%%and $\xi_i^{i^*,d_i}+\xi_{i^*}^{i,d_i}\in\gt g$ for all $i\ne i^*$,
i.e., $\esi(i,i^*,0)=1$.
Thus $\alpha(\widetilde{\gt g}_{e,1})=0$ and, hence, 
$(\gt g_e)_\alpha=\gt g\cap(\widetilde{\gt g}_e)_\alpha$.

Suppose $i\ne i^*$. Then 
$\xi_i^{i,s}+\xi_{i^*}^{i^*,s}\in\gt g$ 
\iff $s$ is odd; and 
$\xi_i^{i^*,s}+\xi_{i^*}^{i,s}\in\gt g$
\iff $s$ is even. Suppose now that 
$i=i^*$. Then 
$\xi_i^{i,s}\in\gt g$ \iff $s$ is odd. 
Summing up, we get 
$\dim(\gt g_e)_\alpha=\frac{1}{2}(\sum\limits_{i=i^*}d_i)+
\sum\limits_{(i,i^*),\ i\ne i^*} d_i=n$.

Next, we show that $\alpha(\gt g_{e,0})=0$.
Since $\sigma(e^s{\cdot}w_i)=
(-1)^se^s{\cdot}\sigma(w_i)=\pm e^s{\cdot} w_i$,
all vectors $\{e^s{\cdot}w_i\}$ are eigenvectors of $\sigma$.
Hence, $\sigma(\xi_i^{j,s})=\pm \xi_i^{j,s}$.
Suppose $i\ne i^*$ and $\sigma(w_i)=w_i$.
Then $\sigma(e^{d_i}{\cdot}w_i)=e^{d_i}w_i$ and, 
since $(e^{d_i}{\cdot}w_i, w_{i^*})\ne 0$, we get
$\sigma(w_{i^*})=-w_{i^*}$, 
$\sigma(e^{d_i}{\cdot}w_{i^*})=-e^{d_i}{\cdot}w_{i^*}$.
Thus $\xi_i^{i^*,d_i}, \xi_{i^*}^{i,d_i}\in\gt g_{e,1}$.
In case $i=i^*$, $d_i$ is odd, 
$\sigma(\xi_i^{i,d_i})=-\xi_i^{i,d_i}$,
and $\xi_i^{i,d_i}\in\gt g_{e,1}$.
Suppose $\varphi=(\sum c_i^{j,s}\xi_i^{j,s})\in\widetilde{\gt g}_{e,1}$.
Then all coefficients $c_i^{i^*,d_i}$ 
of $\varphi$ equal zero. In particular, 
$\alpha(\varphi)=0$. Thus $\alpha(\gt g_{e,0})=0$ and
indeed $\alpha$ is a point of $\gt g_{e,1}^*$. 
Finally, notice that
$\dim(\gt g_{e,1})_\alpha\le\dim(\gt g_e)_\alpha=n$.  
Hence $\dim(\gt g_{e,1})_\alpha=n$, and we are done.
\end{proof}

\noindent
In Section~\ref{ohne}, we show that most of the isotropy representations associated 
with other involutions of $\spv$ do not have GNIB.

                %%%%
%%%%%%%%%%%%%%%%
%%%%%%%%%%%%%%%%   Section 5
%%%%%%%%%%%%%%%%
                %%%%

\section{The isotropy representations of $(\gt{so}_{2n},\gt{gl}_n)$}
%%%symmetric pairs associated with $\sov$}
\label{one-more-2}
\setcounter{equation}{0}

\noindent
In this section, $\dim V=2n$, $\gt g=\gt{so}(V)$, and 
$(\,,\,)$ is a $\g$-invariant symmetric form on $V$.
Let $\sigma$ be an involution of $\g$ such that $\g_0\simeq
\gln$. This can explicitly be described as follows. 
%%%and $\gt g_0=\gt g^\sigma=\gt{gl}_n$. 
Let $V=V_+\oplus V_-$ be a Lagrangian decomposition of $V$.
Then $G_0$ can be taken as the subgroup
of $G=SO(V)$ preserving this decomposition. 
Here $G_0\simeq GL(V_+)$, $V_-\simeq (V_+)^*$ as $G_0$-module, and
the $G_0$-module $\g_1$ is isomorphic to 
$\wedge^2 V_+\oplus (\wedge^2 V_+)^*$.

Keep the notation introduced in Section~\ref{spv}.
In particular, $V_i$, $i=1,\ldots,m$, are the Jordan spaces of $e\in\fN(\glv)$,
$\dim V_i=d_i+1$, and $w_i\in V_i$ is a cyclic vector.

\begin{lm}   \label{restr3}
Let $e\in\fN(\g_1)$.
Then the cyclic vectors $\{w_i\}$ and hence the spaces $\{V_i\}$ 
can be chosen such that the following properties are satisfied:
\begin{itemize}
\item[\sf (i)] \ there is an involution $i\mapsto i^*$ on the set
$\{1,\dots, m\}$ such that 
\begin{itemize}
\item[$\circ$] \ $i^*\ne i$ for each $i$;
\item[$\circ$] \ $d_i=d_{i^*}$; 
\item[$\circ$] \ $(V_i, V_j)=0$ if $i\ne j^*$. In particular, $(V_i,V_i)=0$.
\end{itemize}
\item[\sf (ii)] \ $\sigma(w_i)=\pm w_i$. More precisely, if $d_i$ is even, then the signs
for $\sigma(w_i)$ and $\sigma(w_{i^*})$ are the same; if $d_i$ is odd, then the signs
are opposite.
\end{itemize}
\end{lm}
\noindent 
The proof is left to the reader (cf. the proof of Lemma~\ref{restr}).
%%Actually, part (i) is a standard property of nilpotent orbits in $\spv$.
%%Then part (ii) says that in the presence of the involution $\sigma$
%%%the cyclic vectors for $e\in\fN(\g_1)$ can be chosen to be $\sigma$-eigenvectors.
%%\\[.7ex]

\begin{thm}   \label{more_gib2}
The representation 
$(GL(V_+) : \wedge^2 V_+\oplus (\wedge^2 V_+)^*)$ 
has GNIB. 
\end{thm}
\begin{proof} In the argument below, we omit 
routine but tedious calculations 
of stabilisers and verifications that 
some functions $\alpha \in (\gt g_e)^*$ actually belong to $\gt g_{e,1}^*$.
All this is similar to computations already presented in Sections~\ref{spv} and \ref{one-more}.
%%%analogous ones done in the previous section.  

Recall that 
$\gt{so}(V)$ is a symmetric subalgebra of $\widetilde{\gt g}:=\glv$.  
We follow the notation similar to that used in the proof of Theorem~\ref{more_gib}.
In particular, $\widetilde{\gt g}=\gt{so}(V)\oplus\widetilde{\gt g}_1$ 
is a $\mathbb Z_2$-grading, and there is again a hierarchy of two involutions. 

Let $e\in\fN(\gt g_1)$. In this case, $\rk (G/G_0)=[n/2]$ and,
by Proposition~\ref{z2-gnib}, our goal is to find an element 
$\ap\in (\g_{e,1})^*$ such that $\dim(\g_{e,1})_\ap=[n/2]$.
Choose the cyclic vectors for $e$ as prescribed by
Lemma~\ref{restr3}. 
We normalise these vectors such that 
$(w_i, e^{d_i}{\cdot}w_{i^*})=-(w_{i^*}, e^{d_i}{\cdot}w_i)=\pm 1$.
Then $\gt g_e$ has a basis 
$\xi_i^{j,d_j-s}+\esi(i,j,s)\xi_{j^*}^{i^*,d_i-s}$, 
where $\esi(i,j,s)=\pm 1$ depending on $i,j$ and $s$; and 
$\xi_i^{j,d_j-s}-\esi(i,j,s)\xi_{j^*}^{i^*,d_i-s}$ form a basis 
for $\widetilde{\gt g}_{e,1}$. 

We argue by induction on $m$. Notice that by Lemma~\ref{restr3} $m$ is even.
\\
$\bullet$ \ Suppose first that $m=2$. 
Then $d_1=d_2$ and $(V_i,V_i)=0$. Abusing notation, we write $\sigma(v)/v$ for the sign in
the formula $\sigma(v)=\pm v$.
By Lemma~\ref{restr3}(ii), we have 
$\sigma(w_1)/w_1=\sigma(w_2)/w_2$ if  $d_1$ is odd, and
$\sigma(w_1)/w_1=-\sigma(w_2)/w_2$ if $d_1$ is even. 
%%Hence $\sigma(w_1)/w_1=-\sigma(e^{d_1}\cdot w_2)/(e^{d_1}\cdot w_2)$. 
The algebra $\g_e$ has a basis
\[
\{ \xi_1^{1,s}+(-1)^{s+1}\xi_2^{2,s} 
\mid s=0,\dots,d_1\}\cup
\{ \xi_i^{i,d_1-s}\mid i=1,2;\, 0\le s\le d_1,\, s \text{ is odd}\} \ .
\]
Here $\sigma(\xi_1^{1,s}+(-1)^{s+1}\xi_2^{2,s})=
(-1)^s(\xi_1^{1,s}+(-1)^{s+1}\xi_2^{2,s})$
and $\sigma(\xi_i^{i,d_1-s})=\xi_i^{i,d_1-s}$. Therefore 
$\dim\g_{e,1}=d_1=[n/2]$. Since $\dim(\g_{e,1})_\ap$ cannot be less than
$\rk(G/G_0)$, we obtain $\dim(\gt g_{e,1})_\ap=[n/2]$ for any $\ap$, as required. 
\\[.8ex]
$\bullet$ \ Assume that $m\ge 4$ and the statement holds for all $m_0<m$. 
In the induction step, we use the following simple fact. 
Suppose there is $\alpha\in\gt g_{e,1}^*$ such that
$\ind((\gt g_{e,0})_{\alpha},(\gt g_{e,1})_{\alpha})=[n/2]$.
Then 
\[
[n/2]\le \ind(\gt g_{e,0},\gt g_{e,1})\le 
\ind((\gt g_{e,0})_{\alpha},(\gt g_{e,1})_{\alpha})=[n/2].
\]
Hence, $\ind(\gt g_{e,0},\gt g_{e,1})=[n/2]$.

Choose an ordering of cyclic spaces such that $d_1\ge d_2\ge\dots\ge d_m$.
Without loss of generality,
we may assume that $i^*=i+1$ if $i$ is odd. 
Then there are four possibilities:

{\bf (1)} \ $d_1$ is odd;

{\bf (2)} \ $d_1$ is even and there is some $k\in\{3,\ldots,m-2\}$ such that $d_k$ is also even;

{\bf (3)} \ $d_1,d_2,d_{m-1},d_m$ are even and all other $d_i$ are odd;

{\bf (4)} \ $d_1$, $d_2$ are even and all other $d_i$ are odd. 

\noindent
Consider all these possibilities in turn. In cases {\bf (1)} and {\bf (2)}
we argue by induction, whereas in cases {\bf (3)} and {\bf (4)} we explicitly 
indicate a generic point in $(\g_{e,1})^*$.

\vskip0.5ex
\noindent
{\bf (1)} \,\, Set $\gt f_1=\gt{so}(V_1\oplus V_2)$ and
$\gt f_2=\gt{so}(V_3\oplus\dots\oplus V_m)$. Then 
$e=e_1+e_2$, where $e_i\in\gt f_i$.
Define $\ap\in (\g_{e,1})^*$ by the formula
$\alpha(\varphi)=c_1^{1,d_1}+c_2^{2,d_2}$,
where $c_1^{1,d_1}$,$ c_2^{2,d_2}$
are coefficients of $\varphi\in{\gt g}_e$. 
Then $(\gt g_e)_\alpha=
(\gt f_1)_{e_1}\oplus(\gt f_2)_{e_2}$. 
By the inductive hypothesis,
$\ind((\gt g_{e,0})_{\alpha},(\gt g_{e,1})_{\alpha})=[n/2]$.

\vskip0.5ex
\noindent
{\bf (2)}\,\, Let $k>2$ be the first (odd) number 
such that $d_k$ is even. We may assume that 
$\sigma(w_1)=w_1$ and $\sigma(w_k)=w_k$, while $\sigma(w_2)=-w_2$ and 
$\sigma(w_{k+1})=-w_{k+1}$.
Define $\beta\in (\g_{e,1})^*$ by the formula
$$
\beta(\varphi)=\left(\sum\limits_{i=3}^{k-1} a_ic_i^{i,d_i}\right)+
b_1(c_1^{k+1,d_k}-c_{k}^{2,d_2})+b_2(c_2^{k,d_k}-c_{k+1}^{1,d_1}),
$$
where $a_i=a_j$ if and only if $i=j^*$ and $b_1\ne\pm b_2$. One 
can show that
$(\gt g_e)_{\beta}=\gt h\oplus(\gt f_2)_{e_2}$, where 
$\gt h$ is a subalgebra of $\gt f_1=\gt{so}(V_1\oplus\dots\oplus V_{k+1})$,
$\gt f_2=\gt{so}(V_{k+2}\oplus\dots\oplus V_m)$, 
$e=e_1+e_2$, $e_1\in\gt h$, and $e_2\in\gt f_2$. 
By the inductive hypothesis 
$\ind((\gt f_2)_{e_2,0},(\gt f_2)_{e_2,1})=
\left[\frac{\rk\gt f_2}{2}\right]$. 

Let $\gt p$ be an $(\gt f_1\oplus\gt f_2)$-invariant
complement of $\gt f_1\oplus\gt f_2$ in $\gt g$.
Then there is 
a $\sigma$-invariant decomposition
$\gt g_e=(\gt f_1\oplus\gt f_2)_e\oplus\gt p_e$.
If $\alpha\in(\gt f_1)_e^*$, 
then $(\gt g_e)_\alpha=
((\gt f_1)_e)_\alpha\oplus(\gt f_2)_e\oplus\Ker\hat\alpha$, 
where $\Ker\hat\alpha\subset\gt p_e$
is the kernel of the symplectic form 
$\hat\alpha$ defined by 
$\hat\alpha(\xi,\eta)=\alpha([\xi,\eta])$. 
Since $\Ker\widehat{\beta}=0$, 
this is also true for generic points 
$\alpha\in(\gt f_1)_e^*$ such that $\alpha(\gt g_{e,0})=0$.
Therefore, we can find a point 
$\alpha\in(\gt f_1)_e^*$ such that $\alpha(\gt g_{e,0})=0$ and 
$(\gt g_e)_\alpha=\gt h\oplus(\gt f_2)_e$, where $\ind(\gt h_0,\gt h_1)=
\left[\frac{\rk\gt f_1}{2}\right]$. 
For that point $\alpha$ we get 
$$
\ind((\gt g_{e,0})_\alpha,(\gt g_{e,1})_\alpha)=
\ind(\gt h_0,\gt h_1)+
\ind((\gt f_2)_{e_2,0},(\gt f_2)_{e_2,1})=
\left[\frac{\rk\gt f_1}{2}\right]+\left[\frac{\rk\gt f_2}{2}\right]=
\left[\frac{\rk\gt g}{2}\right].
$$ 

\vskip0.5ex
\noindent
{\bf (3)} \& {\bf(4)}\,\, 
We may assume that $\sigma(w_1)=w_1$ and 
$(w_1,e^{d_1}{\cdot} w_2)=1$.
In case {\bf (3)}, we also assume that 
$(w_{m-1},e^{d_m}{\cdot} w_m)=1$ and $\sigma(w_m)=-w_m$. 
Set $t=m{-}2$ in case {\bf(3)}  and  $t=m$ in case~{\bf(4)}.
Take a point $\alpha\in(\g_{e,1})^*$ such that $\alpha(\varphi)=
b(c_1^{1,d_1-1}+c_2^{2,d_2-1})+\sum\limits_{i=3}^{t} a_ic_i^{i,d_i}$, 
where $a_i=a_j$ if and only if $i=j^*$ and
each $a_i\ne b$. For each $i$ odd, we set 
$\gt h_i:=\gt{so}(V_i\oplus V_{i+1})\cap\gt g_e$. Then there exist
numbers $\esi(1,i),\esi(2,i)\in\{+1,-1\}$, depending on $i$, such that
\[
(\gt g_e)_\alpha= 
\left(\bigoplus\limits_{i \text{ odd}}\gt h_i\right)
\oplus
\left<\xi_1^{i,d_i}+\varepsilon(1,i)\xi_{i^*}^{2,d_1},
      \xi_2^{i,d_i}+\varepsilon(2,i)\xi_{i^*}^{1,d_1}\,|\,
      i=3,4,\dots,m-1,m\right>.
\]
The second summand, denoted by $\gt a$, is a 
commutative ideal of $(\g_e)_\ap$.
%%%$\esi(1,i),\esi(2,i)=\pm 1$, depending on $i$. 
Since $\sigma(w_1)=w_1$ and $\sigma(w_2)=-w_2$,
one of the vectors 
$\xi_1^{i,d_i}+\varepsilon(1,i)\xi_{i^*}^{2,d_1}$,
$\xi_2^{i,d_i}+\varepsilon(2,i)\xi_{i^*}^{1,d_1}$
lies in $\gt g_0$ and another in $\gt g_1$ for each pair $\{i,i^*\}\ne \{1,2\}$ .

Each $\gt h_i$ has a Levi decomposition 
$\gt h_i=\gt l_i\oplus\gt n_i$, where 
$\gt l_i:=\gt h_i\cap\gt{gl}(\bbk w_i\oplus\bbk w_{i^*})$ is reductive
and $\gt n_i$ is the nilpotent radical.
If $d_i$ is even, then $\gt l_i\cong\gt{so}_2$;
and if $d_i$ is odd, then $\gt l_i\cong\gt{sl}_2$.
In any case, $\gt l_i\subset\gt g_0$.  
Moreover, we have $[\gt n_i,\gt a]=0$ for each odd $i$, and 
$[\gt l_i,\xi_1^{j,d_j}+\varepsilon(1,j)\xi_{j^*}^{2,d_1}]=
[\gt l_i,\xi_2^{j,d_j}+\varepsilon(2,j)\xi_{j^*}^{1,d_1}]=0$
if $i\ne 1,j,j^*$.

Set $\gt r:=(\gt g_e)_\alpha$.
We claim that $\ind(\gt r_0,\gt r_1)=[n/2]$.
Define $\beta\in\gt r^*$ by the following rule: 
\\
If $x \in \underset{i}{\bigoplus}\gt h_i$, then $\beta(x)=0$;
if $x$ is one of the vectors 
$\xi_1^{i,d_i}+\varepsilon(1,i)\xi_{i^*}^{2,d_1},
      \xi_2^{i,d_i}+\varepsilon(2,i)\xi_{i^*}^{1,d_1}$,
then $\beta(x)=1$ if $x\in\gt r_1$ and
$\beta(x)=0$ if $x\in\gt r_0$.
In particular $\beta(\gt r_0)=0$, i.e., $\beta\in (\gt r_1)^*$. 
%%$\eta=\xi_1^{i,d_i}+\varepsilon(1,i)\xi_{i^*}^{2,d_1},
%%      \xi_2^{i,d_i}+\varepsilon(2,i)\xi_{i^*}^{1,d_1}$. 
Then 
\[
(\gt r_1)_\beta=\bigoplus_{i\text{ odd}}%%%%%%%^{m-1}
(\gt h_i\cap\gt g_{e,1})\oplus\bbk\eta_0,
\]
where $\eta_0=0$ in case {\bf (4)}, and 
$\eta_0=\xi_1^{m,d_m}-\xi_{m-1}^{2,d_1}+
      \xi_2^{m-1,d_i}-\xi_{m}^{1,d_1}$ in case {\bf (3)}.
For each $\gt h_i$, we have 
$\dim(\gt h_i\cap\gt g_{e,1})=\left[\frac{d_i+1}{2}\right]$.
Therefore, in both cases {\bf (3)} and {\bf (4)} we obtain
$\dim(\gt r_1)_\beta=[n/2]$, as required.
\end{proof}

                %%%%
%%%%%%%%%%%%%%%%
%%%%%%%%%%%%%%%%   Section 6
%%%%%%%%%%%%%%%%
                %%%%

\section{Isotropy representations without GNIB}
\label{ohne}
\setcounter{equation}{0}

\noindent
Here we describe a method for finding nilpotent orbit
in isotropy representations
without equality in \re{vinb=}. There is an obvious method of constructing
${\mathbb Z}_2$-graded Lie algebras: take any $\mathbb Z$-grading and then glue it
modulo 2. This will be applied in the following form.
Given $e\in\fN(\g)$, take an $\tri$-triple
containing $e$, say $\{e,h,f\}$. Consider the $\mathbb Z$-grading of
$\g$ that is determined by $h$:
\[
    \g= \bigoplus_{i\in {\mathbb Z}}\g(i), \text{ where } \
     \g(i)=\{x\in \g\mid [h,x]=ix\} \ .
\]
Here $e\in\g(2)$. Suppose that $e$ is {\it even\/}, i.e., $\g(i)=0$ if
$i$ is odd. Gluing modulo 2 means that we define
$\g_0=\bigoplus_{i\in {\mathbb Z}}\g(4i)$ and
$\g_1=\bigoplus_{i\in {\mathbb Z}}\g(4i+2)$.
Then $e\in \g_1$ and it is sometimes possible to prove that, for this
nilpotent element, Equality~\re{vinb=} does not hold.

Our point of departure is an even nilpotent element $e$ of height 4
(the latter means that $(\ad e)^5=0$).
Then the corresponding $\mathbb Z$-grading is
$\g= \bigoplus_{i=-2}^2 \g(2i)$.
The centraliser of $e$ lies in the non-negative part of this grading, i.e.,
$\g_e=\g(0)_e\oplus \g(2)_e \oplus \g(4)$.  Therefore,
\[
  \g_{e,0}=\g(0)_e\oplus \g(4) \quad \text{and} \quad \g_{e,1}=\g(2)_e \ .
\]
Here $\dim\g(2)_e=\dim\g(2)-\dim\g(4)$ and $\dim\g(0)_e=\dim\g(0)-\dim\g(2)$
and hence 
\[
  \dim\g_1-\dim \g_0=2\dim\g(2)- 2\dim\g(4)-\dim\g(0)=
  \dim\g(2)_e-\dim\g(0)_e-\dim\g(4) \ .
\]
We wish to compare $\ind(\g_0,\g_1)$ and $\ind (\g_{e,0}, \g_{e,1})$.
Let $S$ denote the identity component of a generic stabiliser for 
$(G_0:\g_1)$. Then
\[
   \ind(\g_0,\g_1)=\dim \g_1-\dim\g_0+\dim S=
  \dim\g(2)_e-\dim\g(0)_e-\dim\g(4)+\dim S \ .
\]
In our situation, $\g(4)$ acts trivially on $\g(2)_e$ and hence on
$\g(2)_e^*$. Hence the action {$(G_{e,0}:(\g_{e,1})^*)$} essentially
reduces to  a {\it reductive\/} group action $(G(0)_e:\g(2)_e^*)$.
Let $S^{\{e\}}$ denote the identity component of a generic stabiliser for the
representation $(G(0)_e: \g(2)_e)$.
Then
\[
\ind (\g_{e,0}, \g_{e,1})=\ind (\g(0)_e, \g(2)_e)=
\dim\g(2)_e-\dim\g(0)_e +\dim S^{\{e\}} \ .
\]
Hence
\begin{equation}  \label{differ}
   \delta:=\ind (\g_{e,0}, \g_{e,1})-\ind(\g_0,\g_1)=\dim\g(4)+\dim S^{\{e\}}-\dim S \ ,
\end{equation}
and, as was shown in Proposition~\ref{z2-gnib}, this quantity is non-negative. 
The stabilisers $S$ are well-known. (Actually, they can be directly read off from
the Satake diagram of the involution in question.) Some work is only needed for
computing $\dim S^{\{e\}}$.

\begin{rem}
The involutions obtained in this way are always inner.
\end{rem}

\noindent
Below, we provide a series of examples covered by the previous 
scheme.

\begin{ex}  \label{twice_min}
%%We describe a class of nilpotent elements such that the quantity
%%in Eq.~\re{differ} is positive.
Suppose $\g$ is a simple Lie algebra such that the highest root is a fundamental weight.
Take the weighted Dynkin diagram of the minimal nilpotent orbit.
Then twice this diagram is again a weighted Dynkin diagram.
This new diagram determines an even nilpotent orbit (element) of height 4.
In this situation, $\dim\g(4)=1$ and $\g_0=\g(0)'\oplus \tri$.
Then straightforward calculations show that $S^{\{e\}}=S$.
%%(The latter can be proved a priori, if the highest 
%%root of $\g$ is fundamental, i.e., for $\g\ne \sln, \spn$.)
Hence the quantity in \re{differ} is equal to 1.
This yields the following list of symmetric pairs without GNIB:
\begin{itemize}
\item[] \  $(\GR{E}{8},\GR{E}{7}\times \GR{A}{1})$; \
$(\GR{E}{7},\GR{D}{6}\times \GR{A}{1})$; \ 
$(\GR{E}{6},\GR{A}{5}\times \GR{A}{1})$; \ 
$(\GR{F}{4},\GR{C}{3}\times \GR{A}{1})$; \ 
$(\GR{G}{2},\GR{A}{1}\times \GR{A}{1})$;
\item[]  \ $(\gt{so}_n, \gt{so}_{n-4}\times\gt{so}_4)$, $n\ge 7$; 
%%%$(\gt{sl}_n, \gt{sl}_{n-2}\times\gt{sl}_2\times T_1)$, $n\ge 4$;
\end{itemize}
\end{ex}
\begin{rem}
If $\g$ is of type $\GR{G}{2}$, then this procedure leads to Example~\ref{contr}.
%%%in Section~\ref{GIB:sect}.
\end{rem}
\begin{ex}  \label{ex-gl}
Let $e$ be a nilpotent element in $\gt{gl}_{3k+l}$ corresponding to
the partition $(3^k,1^l)$.
%% with $n\ge 2$. 
Then $e$ is even and of height 4, and the related
symmetric pair is $(\gt{gl}_{3k+l}, \gt{gl}_{2k}\times \gt{gl}_{k+l})$.
We have the following data for the dimension of
graded pieces for the $\mathbb Z$-grading:
\begin{center}
\begin{tabular}{c|ccc}
$i$             &  0             &  2        & 4  \\ \hline
$\dim \g(i)$    & $2k^2+(k+l)^2$ & $2k(k+l)$ & $k^2$ \\
$\dim \g(i)_e$  & $k^2+l^2$      & $k^2+2kl$ & $k^2$
\end{tabular}
\end{center}
To compute $S^{\{e\}}$, we notice that 
$\g(0)_e \simeq \gt{gl}_k\times \gt{gl}_l=\gt{gl}(V_1)\times\gt{gl}(V_2)$
and the $\g(0)_e$-module $\g(2)_e$ is isomorphic to 
$(V_1\otimes V_1^*) \oplus (V_1\otimes V_2) \oplus (V_1\otimes V_2)^*$.
Suppose $k\ge l$. Then $S^{\{e\}}=T_1\times GL_{l-k}$,
where $T_{j}$ stands for a $j$-dimensional torus. 
%%It follows that if  $k\ge 2$.
The group $S$ is isomorphic to $T_{2k}\times GL_{l-k}$, 
Hence the quantity $\delta$ in \re{differ} is equal to $k^2-2k+1$, which is positive
for $k\ge 2$. The same type of argument shows that $\delta=1$ if 
$k=2$ and $l=1$.  In particular, this means that 
%%the isotropy representation of
the symmetric pair $({\frak gl}_n,
{\frak gl}_4\times {\frak gl}_{n-4})$ has no GNIB for any $n\ge 7$.
\end{ex}
\begin{ex}  \label{ex-sp}
Now $\g=\spv$.
Let $e$ be a nilpotent element in $\gt{sp}_{6k+2l}$ corresponding to
the partition $(3^{2k},1^{2l})$. 
Then $e$ is even and of height 4, and the related
symmetric pair is $(\gt{sp}_{6k+2l}, \gt{sp}_{4k}\times \gt{sp}_{2k+2l})$.
Here $\g(0)_e=\gt{sp}_{2k}\times \gt{sp}_{2l}=\gt{sp}(V_1)\times \gt{sp}(V_2)$ 
and $\g(2)_e$
is isomorphic to $\wedge^2 V_1\oplus (V_1\otimes V_2)$. 
%%where $\dim V_1=2$, $\dim V_2=2l$,
%%and $\odin$ stands for the trivial one-dimensional module.
If $k\le l$, then  $S^{\{e\}}= \left\{\begin{array}{rc}
Sp_{2(l-k)}, & \text{if } \ k\ge 3 \\
SL_2 \times Sp_{2(l-k)},    & \text{if } \ k\le 2 \end{array}\right.$. 
Also, $S=(SL_2)^{2k}\times Sp_{2(l-k)}$.
%%Therefore
Since $\dim\g(4)=k(2k+1)$, we see that $\delta$
%%%the quantity in \re{differ} 
is positive for $k\ge 2$. Similarly to the previous example, one also verifies
that $\delta$ is positive for $k=2, l=1$. 
In particular,  
%%the isotropy representation of
the symmetric pair $({\frak sp}_{2n},
{\frak sp}_8\times {\frak sp}_{2n-8})$ has no GNIB for any $n\ge 7$.
\end{ex}
\begin{ex}  \label{ex-e7}
Let $\g$ be of type $\GR{E}{7}$ and $e$ a nilpotent element with weighted Dynkin
diagram 
\begin{tabular}{@{}c@{}}
0-0-0-\lower3.1ex\vbox{\hbox{0\rule{0ex}{2.4ex}}
\hbox{\hspace{0.5ex}\rule{.1ex}{1ex}\rule{0ex}{1.2ex}}\hbox{2\strut}}-0-0
\end{tabular}. Then $\dim\g(0)=49$, $\dim\g(2)=35$, and $\dim\g(4)=7$.
Therefore $\dim\g_0=63$ and $\dim\g_1=70$. Hence the related involution is
of maximal rank. 
(Here $\g_0\simeq \gt{gl}(V)$ with $\dim V=7$ and $\g_1\simeq \wedge^4 V$). 
Therefore the group $S$ is trivial. This already means that
$\dim\g(4)+\dim S^{\{e\}}-\dim S\ge 7$.
\end{ex}
\noindent

\begin{ex}  \label{ex-e8}
Let $\g$ be of type $\GR{E}{8}$ and $e$ a nilpotent element with weighted Dynkin
diagram 
\begin{tabular}{@{}c@{}}
0-0-0-0-\lower3.1ex\vbox{\hbox{0\rule{0ex}{2.4ex}}
\hbox{\hspace{0.5ex}\rule{.1ex}{1ex}\rule{0ex}{1.2ex}}\hbox{0\strut}}-0-2
\end{tabular}. Then $\dim\g(0)=92$, $\dim\g(2)=64$, and $\dim\g(4)=14$.
Therefore $\dim\g_0=120$ and $\dim\g_1=128$. Hence the related involution is
of maximal rank. 
(Here $\g_0\simeq \gt{so}_{16}$ and $\g_1$ is a half-spinor representation of
$\gt{so}_{16}$). 
Therefore the group $S$ is trivial. This already means that
$\dim\g(4)+\dim S^{\{e\}}-\dim S\ge 14$.
\end{ex}

\noindent
One can find more isotropy representations without GNIB using the above examples
and the slice method. The following assertion readily follows
from Theorems~\ref{gnib+gib} and \ref{oboz}.
\\[.7ex]
\hbox to \textwidth{\refstepcounter{ex} (\theex) \quad  \label{claim}
\parbox[t]{422pt}{{\sl Suppose $(G:V)$ is observable and $(L:W)$ is a slice representation
of $(G:V)$. If $(L:W)$ has no GNIB, then so does $(G:V)$.}}
}

\begin{ex}  \label{ex-e6}
Consider the symmetric pair $(\GR{E}{6},{\frak sp}_8)$. The corresponding
$\mathbb Z_2$-grading $\g=\g_0\oplus\g_1$ is of maximal rank
and the $Sp_8$-module $\g_1$ is isomorphic to the 4-th fundamental representation.
Since the rank is maximal, for any Levi subalgebra $\el\subset \GR{E}{6}$,
there is $x\in \g_1$ whose centraliser in $\g$ is conjugate to $\el$. 
The induced $\mathbb Z_2$-grading of $\z(x)\simeq\el$ is also of maximal rank.
%%It follows that the isotropy representation of the $\mathbb Z_2$-grading of 
%%maximal rank of $\el$ occurs as a slice representation of $(G_0:\g_1)$.
In particular, taking $x$ such that the semisimple part of $\z(x)$ is of type
%%%for the Levi subalgebra whose semisimple part is of type
$\GR{D}{4}$, we obtain, up to the centre of $\z(x)$, the symmetric
pair $({\frak so}_8,{\frak so}_4\times {\frak so}_4)$, which has no GNIB.
Hence the symmetric pair $(\GR{E}{6},{\frak sp}_8)$
has no GNIB, too. \\
The similar argument also works for the involutions in Examples~\ref{ex-e7},\,\ref{ex-e8}.
\end{ex}

\begin{thm}  \label{more}
For the following symmetric pairs, the isotropy representation has no GNIB:
\begin{itemize}
\item[\sf(i)] \ $(\gt{gl}_{n}, \gt{gl}_{m}\times \gt{gl}_{n-m})$ with $4\le m\le n-m$;
\item[\sf(ii)] \ $(\gt{so}_{n}, \gt{so}_{m}\times \gt{so}_{n-m})$ with $4\le m\le n-m$;
\item[\sf(iii)] \ $(\gt{sp}_{2n}, \gt{sp}_{2m}\times \gt{sp}_{2n-2m})$ with 
$4\le m\le n-m$.
%%\item[\sf(iv)] \ 
\end{itemize}
The cases with $m=3$ and $n-m=4$ also yield the isotropy 
representations without GNIB.
\end{thm}\begin{proof}
(i) It is easily seen that $(G_0:\g_1)$ has a slice representation
which is isomorphic to the isotropy representation of the pair
$({\frak gl}_{n-2},{\frak gl}_{m-1}\times {\frak gl}_{n-m+1})$.
Iterating this procedure yields 
the isotropy representation of the pair 
$({\frak gl}_{n-2m+8},{\frak gl}_{4}\times {\frak gl}_{n-2m+4})$. The latter has no GNIB
by Example~\ref{ex-gl}. Then one applies assertion~\ref{claim}.
\\[.6ex]
(ii),(iii). Here one uses the similar reductions, with `$\gt{so}$' and 
`$\gt{sp}$' in place of `$\gt{gl}$', and Examples~\ref{twice_min} and \ref{ex-sp}. 
\end{proof}%
\noindent
Making use of a direct computation, we strengthen the assertion 
of Theorem~\ref{more}(i).

\begin{thm}            \label{rk3}
The following symmetric pairs $(\g,\g_0)$ have no GNIB:
\begin{itemize}
\item[\sf (i)] \ 
$(\gt{gl}_{n}, \gt{gl}_3\oplus\gt{gl}_{n-3})$, $n\ge 7$,
\item[\sf (ii)] \ $(\gt{so}_{n}, \gt{so}_3\oplus\gt{so}_{n-3})$, $n\ge 7$,
\item[\sf (iii)] \ $(\gt{sp}_{2n}, \gt{sp}_6\oplus\gt{sp}_{2n-6})$, $n\ge 7$.
\end{itemize} 
\end{thm}
\begin{proof}  
For all these symmetric pairs, we have $\rk(G/G_0)=3$, and the case of $n=7$ is covered
by Theorem~\ref{more}. We show that for $n > 7$ there is a reduction to $n=7$.

(i) Let $\gt h=\gt{gl}_7\subset\gt{gl}_{n}$ 
be a regular $\sigma$-invariant subalgebra 
such that $\gt h^{\sigma}=\gt{gl}_3\oplus\gt{gl}_4
\subset \gt{gl}_3\oplus\gt{gl}_{n-3}$. 
By Example~\ref{ex-gl}, the nilpotent $H$-orbit with partition $(3,3,1)$
meets $\h_1$ and yields an $H_0$-orbit without GNIB. Let $e\in\h_{1}$ be an element in 
this orbit.
Using the embedding $\h_1\subset \g_1$, we may regard $e$ as element of $\g_1$.
Then the corresponding partition is $(3,3,1^{n-6})$.
We are going to prove that $\ind(\gt g_{e,0},\gt g_{e,1})\ge \ind(\gt h_{e,0},\gt h_{e,1})$.
%%According to Example~\ref{ex-gl}, $\ind(\gt h_{e,0},\gt h_{e,1})>3$. 
%%We claim that $\ind(\gt g_{e,0},\gt g_{e,1})>3$. 
An explicit model of $e$ is as follows. Let $w_1,w_2,\ldots,w_{n-4}$ be cyclic vectors
for $e$, where $e^{3}{\cdot}w_i=0$ for $i=1,2$ and $e{\cdot}w_i=0$ for $i\ge 3$.
Let $\bbk^n=\bbk^3\oplus\bbk^{n-3}$ be the $\g_0$-stable 
decomposition corresponding to $\sigma$. Then we assume that $w_3\in\bbk^3$ and 
and all other $w_i$'s lie in $\bbk^{n-3}$. 
(Hence $\bbk^3=\langle e{\cdot}w_1, e{\cdot}w_2, w_3 \rangle$.)
Now, all information for $e$ can be presented rather explicitly.
We have
%%One can calculate dimensions of various stabilisers:
$\dim\gt g_{e,0}=(n{-}4)^2+5$,
$\dim\gt h_{e,0}=9$, $\dim\gt h_{e,1}=8$.
More precisely,
$\gt g_{e,0}=\gt h_{e,0}\oplus\gt{gl}_{n-7}\oplus\gt a$, 
where $[\gt{gl}_{n-7},\gt h]=0$ and 
\[
\gt a=\langle \xi_1^{t,0},\xi_2^{t,0},\xi_t^{1,2},\xi_t^{2,2}\mid 4\le t\le n-4\rangle .
\]
This means in particular that $\ah=\{0\}$ if $n=7$. Next,
\[
\gt h_{e,1}=\left<\xi_i^{j,1},\xi_i^{3,0},\xi_3^{i,2}\mid i,j\in\{1,2\} \right>.
\]
It is easily seen that 
$[\gt a,\gt h_{e,1}]=0$.
For instance, $[\xi_t^{1,2},\xi_1^{i,1}]=-\xi_t^{i,3}=0$, since 
$e^3\cdot w_i=0$.
Because also
$[\gt{gl}_{n-4},\gt h_{e,1}]=0$, 
we get 
$[\gt g_{e,0},\gt h_{e,1}]=[\gt h_{e,0},\gt h_{e,1}]\subset\gt h_{e,1}$. 
It remains to observe that for each $\alpha\in(\gt g_{e,1})^*$ 
we have 
\[
\dim(\gt g_{e,1})_\alpha\ge
\dim (\gt h_{e,1})_\alpha =
\dim (\gt h_{e,1})_{\widetilde\alpha}\ge 4 > \rk(G/G_0),
\]
where $\widetilde\alpha\in(\gt h_{e,1})^*$ is the restriction
of $\alpha$. 

(ii) \ The previous argument goes through {\sl mutatis mutandis\/} 
in the orthogonal case. We just consider the nilpotent element in $\gt{so}_7$
with partition $(3,3,1)$.  Using the natural embedding 
$\gt{so}_7 \subset\gt{so}_n$, $n>7$, we obtain the nilpotent orbit with partition
$(3,3,1^{n-6})$.

(iii) \,This case is similar to part (i) but in a different fashion.
Here we start with a nilpotent element in $\gt{sp}_{14}$ with partition 
$(3,3,3,3,1,1)$, and then embed $\gt{sp}_{14}$ in $\gt{sp}_{2n}$, $n>7$.
\end{proof}
\noindent
%%It is not claimed in the theorem that
%%$(\gt{gl}_6,\gt{gl}_3\times\gt{gl}_3)$ has no GNIB. Actually, one can
%%verify that this symmetric pair has GNIB.

                %%%%
%%%%%%%%%%%%%%%%
%%%%%%%%%%%%%%%%   Section 7
%%%%%%%%%%%%%%%%
                %%%%

\section{More affirmative results and open problems}
\label{optimistic_end}
\setcounter{equation}{0}

\noindent
We conclude with some more examples of isotropy representations having GNIB
and state several questions.

\begin{prop}   \label{rank1}
If\/ $(G,G_0)$ is a symmetric pair of rank 1, then $(G_0:\g_1)$ has GNIB.
\end{prop}\begin{proof}
The symmetric pairs of rank 1 are the following:
\begin{itemize}
\item[]  \ $(\gt{gl}_n, \gt{gl}_{n-1}\times \gt{gl}_{1})$, \ 
$(\gt{sp}_{2n}, \gt{sp}_{2n-2}\times\gt{sp}_2)$, \ 
$(\gt{so}_n, \gt{so}_{n-1})$, \ 
$(\GR{F}{4},{\frak so}_9)$. 
\end{itemize}
The number of nonzero nilpotent $G_0$-orbits in $\fN(\g_1)$
equals $3,\,2,\,1,\,2$, respectively. By Proposition~\ref{max_orb}, the orbit 
of maximal dimension is always "good", so that it remains to test the minimal orbit(s).
This is done by hand. 

We give some details for the last case.
Here the isotropy representation is the spinor (16-dimensional) representation of 
$\g_0={\frak so}_9$.
The weights are $\frac{1}{2}(\pm\esi_1\pm\esi_2\pm\esi_3\pm\esi_4)$. 
For simplicity, weights will be represented by the set of 4 signs.
For instance, the lowest weight is $(----)$. 
Let $v\in \g_1$ be a lowest weight vector.
Then $\dim\g_0{\cdot}v=11$ and $(\g_0)_v$ is a semi-direct product of $\gt{sl}_4$ 
and a nilpotent radical. As $\gt{sl}_4$-module, the 5-dimensional
space $\g_1/\g_0{\cdot}v$ can be identified with the subspace $W$ of $\g_1$ whose
weights are $(-+++),\,(+-++),\,(++-+),\,(+++-),\,(++++)$. Hence $W$ is the sum
of the trivial and 4-dimensional $\gt{sl}_4$-modules. This shows that 
already $SL_4$, the reductive part of $(G_0)_v$, has an orbit of codimension one
in $\g_1/\g_0{\cdot}v$.
\end{proof}
\begin{ex} 
The symmetric pair $(\GR{E}{6},\GR{F}{4})$ has rank two.
However, its isotropy representation 
has only two nonzero nilpotent orbits.
Here again one can easily check that the minimal orbit $\mathcal O_{min}$ satisfies
GNIB-condition, i.e., Equality~\re{vinb=} is satisfied for $v\in \mathcal O_{min}$. 
Hence this isotropy representation has GNIB.
\end{ex}
\begin{rem}      \label{rk2-gl}
Using explicit description of nilpotent $G_0$-orbits, one can honestly verify that
the symmetric pairs $(\gt{gl}_{n}, \gt{gl}_2\oplus\gt{gl}_{n-2})$ and
$(\gt{gl}_{6}, \gt{gl}_3\oplus\gt{gl}_{3})$ have GNIB.
Together with results of Sections~\ref{spv} and \ref{ohne},
this completes the problem of classifying the isotropy representations
of $\gt{gl}_n$ with and without GNIB.
\end{rem}
\begin{rem}      \label{rk2-so}
It is also not hard to verify that the pairs 
$(\gt{so}_{n}, \gt{so}_{2}\times \gt{so}_{n-2})$  and
$(\gt{sp}_{2n}, \gt{sp}_{4}\times \gt{sp}_{2n-4})$ have GNIB. 
Furthermore, both pairs $(\gt{so}_{6}, \gt{so}_{3}\times \gt{so}_{3})$
and $(\gt{sp}_{12}, \gt{sp}_{6}\times \gt{sp}_{6})$
have GNIB. Together with results of Sections~\ref{one-more}--\ref{ohne},
this completes the problem of classifying the isotropy representations
of $\gt{so}_n$ and $\gt{sp}_{2n}$ with and without GNIB.
\end{rem}
Taking into account all symmetric pairs considered so far,
one may notice that there remain only two unmentioned symmetric pairs:
$(\GR{E}{6},{\frak so}_{10}\times \te_1)$ and $(\GR{E}{7},\GR{E}{6}\times \te_1)$.
Their ranks are 2 and 3, respectively.
It is likely that the first of them has GNIB, but we have no assumption
for the second case.
We hope to consider these remaining cases in a subsequent article.

\noindent
There are many interesting open questions on GIB and GNIB. Here are some of them.
\begin{itemize}
\item[(Q1)] \  We have shown in Corollary~\ref{GIB:big} that sufficiently large
reducible representations have GIB. However, no a priori results is known
for irreducible representations of simple algebraic groups.
We conjecture that for any semisimple $G$ there are finitely many 
irreducible representations without GNIB.
%%E.g., is it true that the number of irreducible representations with
%%(or without) GIB is finite?
\item[(Q2)] \ Let $V$ be a simple $G$-module and $v\in V$ a highest weight vector.
Is it true that Equality~\re{vinb=} holds for $v$?
\item[(Q3)] \ Suppose $G$ has a dense orbit in $V$, i.e.,
$\bbk(V)^G=\bbk$. Is it true that  $V$ has GNIB ?
\item[(Q4)] \ Let $V$ be a spherical $G$-module. Is it true that $V$ has GNIB ?
(It is a special case of (Q3).)
\end{itemize}

\noindent
In connection with the last question, we mention that most of spherical 
modules are obtained by the following construction.
Let $\p\subset \g$ be a parabolic subalgebra whose nilpotent
radical, $\p^u$, is Abelian. Let $\p=\el\oplus \p^u$ be a Levi decomposition.
Then $\p^u$ is a spherical $L$-module.
Using the theory developed in \cite{canad}, one can prove that $\p^u$ has GNIB.
The point here is that, for any $v\in\p^u$, already the reductive part of
$L_v$ has an open orbit in $\p^u/\el{\cdot}v$.

Finally, we recall that most of the observable representations of reductive groups
are associated with automorphisms of finite order of simple Lie algebras,
i.e., the corresponding linear groups are $\Theta$-groups in the sense
of Vinberg \cite{vi}. This is a generalisation of the situation considered in this paper.
It is therefore natural to investigate
when these $\Theta$-representations have GNIB.

\end{document}